\documentclass[10pt]{amsart}
\usepackage{amsmath, amsfonts,amssymb}
\usepackage{amsfonts,mathrsfs, color, amsthm}
\usepackage{mathtools}
\usepackage{dsfont} 
\usepackage{adjustbox}
\usepackage{enumerate}
\usepackage[normalem]{ulem}
\usepackage{amsrefs}
\usepackage{tikz}
\usepackage{comment}
\usepackage{ragged2e}
\usetikzlibrary{positioning,patterns,calc}
\DeclareMathAlphabet{\mathpzc}{OT1}{pzc}{m}{it}
\usepackage{float}
\usepackage{graphicx}
\usetikzlibrary{positioning,patterns,calc}

\newcommand{\real}{\mathbb{R}} 
\newcommand{\naturals}{\mathbb{N}} 
\newcommand{\nball}{B^N} 

\newcommand{\divergent}{\textup{div}} 
\newcommand{\graph}{\textup{graph}}
\newcommand{\area}{\textup{area}}
\newcommand{\volume}{\textup{vol}}
\newcommand{\vol}{\textup{vol}}

\newcommand{\support}{\textup{spt}}

\newtheorem{theorem}{Theorem}[section]
\newtheorem{lemma}[theorem]{Lemma}
\newtheorem{remark}[theorem]{Remark}
\newtheorem{proposition}[theorem]{Proposition}

\newtheorem{definition}[theorem]{Definition}
\newtheorem*{theorem*}{Theorem}
\newtheorem{innercustomthm}{Theorem}
\newenvironment{maintheorem}[1]
  {\renewcommand\theinnercustomthm{#1}\innercustomthm}
  {\endinnercustomthm}
\newtheorem{claim}{\texttt{Claim}}
\makeatletter
\@addtoreset{claim}{theorem}
\makeatother
\newtheorem{subclaim}{\texttt{Claim}}[claim]

\usepackage{etoolbox}
\makeatletter
\patchcmd{\@maketitle}
  {\ifx\@empty\@dedicatory}
  {\ifx\@empty\@date \else {\vskip3ex \centering\footnotesize\@date\par\vskip1ex}\fi
   \ifx\@empty\@dedicatory}
  {}{}
\patchcmd{\@adminfootnotes}
  {\ifx\@empty\@date\else \@footnotetext{\@setdate}\fi}
  {}{}{}
\makeatother
\makeatletter
\def\specialsection{\@startsection{section}{1}%
  \z@{\linespacing\@plus\linespacing}{.5\linespacing}%
  {\normalfont}}
\def\section{\@startsection{section}{1}%
  \z@{.7\linespacing\@plus\linespacing}{.5\linespacing}%
  {\normalfont\scshape}}
\makeatother

\title[Compactness of Free Boundary CMC Surfaces]{Compactness of the Space of Free Boundary CMC Surfaces with Bounded Topology}
\author{Nicolau S. Aiex and Han Hong}
\date{\today}
\address{Department of Mathematics, National Taiwan Normal University, Taipei, Taiwan}
\email{nsarquis@math.ntnu.edu.tw}
\address{Yau Mathematical Sciences Center\\ Jing Zhai\\ Tsinghua University\\ Hai Dian District\\ Beijing China 100084}
\email{hh0927@tsinghua.edu.cn}

\begin{document}

\begin{abstract}
We prove that the space of free boundary CMC surfaces of bounded topology, bounded area and bounded boundary length is compact in the $C^k$ graphical sense away from a finite set of points.
This is a CMC version of a result for minimal surfaces by Fraser-Li \cite{fraser.a-li.m2014}. 
\end{abstract}

\maketitle

\section{Introduction}

A Constant Mean Curvature (CMC) surface $\Sigma$ is a critical point of the area functional with respect to variations that preserve enclosed volume.
As a consequence of the first variation of area, the scalar mean curvature has to be constant for a given choice of normal direction.
If $\Sigma$ is an immersed surface with boundary immersed in $N$ with boundary and $\partial \Sigma \subset \partial N$, we say the surface is free boundary CMC if it is a critical point with respect to variations that are in addition tangent along $\partial N$.

Although most results about minimal surfaces have equivalent versions for CMC surfaces, there are many distinctions between their behaviour.
For example, when the mean curvature is non-zero, the mean curvature vector defines a trivialization of the normal bundle.
That is, every CMC surface in a $3$-manifold is $2$-sided.
On the other hand, CMC surfaces may have tangential self-touching points as long as the mean curvature vector points at opposite directions on those points.

The goal of this article is to prove the CMC equivalent to Fraser-Li's result \cite{fraser.a-li.m2014}.
We prove:

\begin{maintheorem}{\ref{compactness fbcmc}}
Let $N$ be a compact $3$-dimensional manifold with boundary.
Suppose $H_{\partial N}> H_0$ with respect to the inward conormal of $N$ along $\partial N$ and let $\Sigma_i$ be a sequence of {connected} embedded free boundary CMC surfaces with mean curvature $H_i$, genus $g_i$ and number of ends $r_i$ satisfying:
\begin{enumerate}[(a)]
\item $|H_i|\leq H_0$;
\item $g_i\leq g_0$;
\item $r_i\leq r_0$;
\item $\textup{area}(\Sigma_i)\leq A_0$ and
\item $\textup{length}(\partial \Sigma_i)\leq L_0$.
\end{enumerate}

Then there exists a smooth properly almost embedded CMC surface $\Sigma\subset N$ and a finite set $\Gamma\subset \Sigma$ such that, up to a subsequence, $\Sigma_i$ converges to $\Sigma$ locally graphically in the $C^k$ topology on compact sets of $N\setminus \Gamma$ for all $k\geq 2$.
Moreover, if $\Sigma$ is minimal then it is properly embedded.

If in addition $\partial N$ satisfies $A^{\partial N}> 0$ with respect to the inward conormal of $\partial N$, then the convergence is $1$-sheeted away from $\Gamma$ whenever $H_\Sigma\neq 0$.
\end{maintheorem}

We note that the theorem holds true when $\partial N=\partial \Sigma=\emptyset$.
In this case one simply ignores all boundary hypotheses.

Let us highlight the main differences that justify the extra hypotheses and the weaker compactness for embedded surfaces.

Firstly, we mention that there is no Steklov eigenvalue estimates for free boundary CMC surfaces.
Consequently, the hypothesis of length bound remains crucial.
Secondly, there is no suitable isoperimetric inequality that allows us to remove the bound on the area.
One could exchange the length bound condition by stability of the surface (see \cite{mendes2018}) but as seen on \cite{aiex-hong2021}, in particular, stable free boundary CMC surfaces have bounded topology so this condition would be topologically restrictive.
Thirdly, even under positive Ricci curvature of the ambient space a sequence of free boundary CMC surfaces may have a neck-pinching phenomenum where the norm of the second fundamental form blows-up.
Naturally the convergence is not smooth along these points where curvature is accumulating.
Finally, the maximum principle for CMC surfaces only applies when their mean curvature vectors point in the same direction.
That is, a sequence of embedded free boundary CMC surfaces may touch tangentially in the limit as long as the limiting surface is not minimal.

Despite the lack of certain useful technical results, the proof of the theorem relies on the same main ideas, each of which is proved using a less optimal technique to compensate for the above limitations.
These are: $L^2$-curvature bounds from the topological bound, improvement to local pointwise curvature estimates, a removable singularity theorem for interior and boundary points and construction of weakly stable free boundary CMC surfaces to study the multiplicity of the limit.

Let us briefly address our approach to each of the above steps.
The integral curvature bounds follow directly from Gauss-Bonnet Theorem, from which the fact the mean curvature is non-zero introduces an extra term of area and the geodesic curvature along the boundary brings in a term of boundary length.
This is essentially the only step in which dimension $2$ is relevant.
The local pointwise curvature estimate comes from a blow-up argument as in \cite{white1987} and Schauder estimates.
It becomes relevant that the blow-up of a CMC surface around a point is a minimal surface in $\real^3$.
Unlike \cite{fraser.a-li.m2014}*{Theorem 4.1} we do not know whether a CMC surface is conformally equivalent to a punctured Riemann surface so we do not have the branch point structure to prove the removable singularity theorem.
This will follow from a blow-up argument to prove that tangent cones are in fact planes (or half-planes) and a local free boundary CMC foliation argument to prove that it is also unique.
These are the same ideas as in \cite{white1987} together with the methods in \cite{ambrozio-carlotto-sharp2018.3} to deal with boundary singularity points.
In fact, most of the calculations are done in the latter reference and only minor adaptations are necessary for the non-minimal case.

The extra condition on the mean curvature of $\partial N$ is only necessary to apply the maximum principle for CMC surfaces and ensure that the interior of the limiting surface is properly immersed in $N$.
Removing this condition would allow for interior tangential touching points between $\partial N$ and the limit surface.

Finally, we discuss the $1$-sheeted convergence part of the statement.
Unlike the minimal case \cite{fraser.a-li.m2014}, we do not require a positive Ricci curvature condition.
The multiplicity analysis in the free boundary case is similar to the closed case in \cite{bourni-sharp-tinaglia}*{\textsection 4}.
We follow the original ideas in the proof of \cite{bourni-sharp-tinaglia}*{Theorem 4.1} but we address a minor mistake found in the proof of \cite{bourni-sharp-tinaglia}*{Claim 4.5} (see also Remark \ref{mistake 1 remark}).

We also remark that our result not only generalizes the Compactness Theorem of \cite{sun.a2020}*{Theorem 1.1} to the free boundary case but also improves it by removing the positive Ricci curvature condition.
In fact we should point out that there is an error in the proof of $1$-sheeted convergence in \cite{sun.a2020}*{Theorem 4.5} which we also fix.
More specifically, there is an mistake when proving that a positive solution to the Jacobi equation can be extended across points of singular convergence.
The application of the maximum principle to prove that such solutions are uniformly bounded is in fact incorrect.
The author does not verify that the first point of touch of the foliation is in fact in the interior and it turns out that it could be at a boundary point.
Moreover, unlike it is claimed in \cite{sun.a2020}*{Theorem 1.1(2)}, it is not clear whether the number of sheets can be controlled when the limiting surface is minimal (see Remark \ref{mistake 2 remark}).
Prior to knowing that the convergence is $1$-sheeted, it is not yet known that there exists a non-trivial Jacobi vector field on the limiting surface.

This article is divided as follows. 
Section $2$ establishes notation, necessary definitions and we prove the geometric version of Schauder estimates.
In section $3$ we prove the local pointwise curvature estimates from uniform integral curvature bounds and the corresponding version for weakly stable CMC surfaces.
The Removable Singularity Theorem is proved in section $4$ and we write the proof for the specific case of a singularity along the boundary.
Section $5$ is dedicated to prove the Compactness Theorem.


\textit{
Acknowledgements: The majority of the work on this article took place when the first author was a postdoctoral fellow at PIMS-University of British Columbia.
The authors would like to thank Professors Jingyi Chen and Ailana Fraser for useful discussions and support.
The authors are also grateful to the authors of \cite{bourni-sharp-tinaglia} for explaining their paper.
N.A. is funded by NSTC grant 112-2811-M-003-002. H.H. is supported by NSFC Grant No. 12401058 and the Talent Fund of Beijing Jiaotong University No. 2024XKRC008.
}


\section{Preliminaries}

In this section we establish notation and preliminary results that will be used throughout this article.

Let $(N^3,\partial N, g)$ be a $3$-dimensional manifold with non-empty boundary and a smooth Riemannian metric $g$.

\begin{definition}\label{almost embedded}
We say that an immersed surface $\Sigma\subset N$ with non-empty boundary is properly almost embedded in an open set $U\subset N$ if $(\Sigma\setminus\partial \Sigma)\cap U$ is properly immersed in $(N\setminus \partial N)\cap U$, $\partial \Sigma\subset \partial N\cap U$ and there exists a set $\mathcal{S} \subset\Sigma\cap U$ such that
\begin{enumerate}[(a)]
\item $\Sigma\setminus\mathcal{S}$ is embedded;
\item for each $p\in\mathcal{S}$ there exists a neighbourhood $V$ of $p$ in $U$ such that $\Sigma\cap V$ is a union of connected components $W_j$, $j=1,\ldots,l_p$, each $W_j$ is embedded, for each $j\neq j'$ we have $W_j$ lying to one side of $W_{j'}$ and $W_j\cap W_{j'}\subset \mathcal{S}\cap V$.
\end{enumerate}
In particular each $W_j$ can be written as a graph over $W_{j'}$ in a sufficiently small neighbourhood of $p\in\mathcal{S}$.
The set $\mathcal{S}$ is called the self-touching set of $\Sigma$.
We say $\Sigma$ is free boundary if in addition $\Sigma$ is orthogonal to $\partial N$.
\end{definition}


The following lemma is a standard application of Schauder estimates in geometry.
It will allow us to improve from $C^{1,\alpha}$ to $C^{2,\alpha}$ graphical convergence of surfaces, given pointwise curvature estimates.
For the sake of simplicity, given a point $y\in\partial\Sigma$ with inward conormal vector $\nu(y)$ we denote $T^+_y\Sigma=\{v\in T_y\Sigma:g^\Sigma(v,\nu(y))\geq 0\}$ and $B^+_r(0)=B_r(0)\cap T^+_y\Sigma$ the half ball of radius $r$ centered at $0$ with respect to the metric $g^\Sigma_y$, as well as the equivalent notation for the ambient manifold $N$.
Observe that when $\Sigma$ is free boundary, the inward conormal of $\partial N$ and $\partial \Sigma$ coincide.

\begin{lemma}\label{schauder estimates}
Let $(N, g^N)$ be a smooth three dimensional compact Riemannian manifold with boundary, $A>0$ and $f\in C^{0,\alpha}(N)$ with $|f|_{0,\alpha}<H_0$.
There exists $r_0(N, A)>0$ and $C_0(N, A, H_0, \alpha)>0$ such that for any $C^2$ two-sided {isometrically} immersed compact surface with free boundary $\varphi:(\Sigma,g^\Sigma)\rightarrow (N,g^N)$ with mean curvature $H_{\Sigma}=f$ and $x_0\in \Sigma$ satisfying
\begin{equation*}
|A^\Sigma(x)|\leq A  \text{ on } B^{\Sigma}_{r_0}(x_0),
\end{equation*}
there {exist} $y_0\in B^{\Sigma}_{r_0}(x_0)$ and a function $u$ defined on $\Omega_0\subset \varphi_*T_{y_0}\Sigma\subset T_{q_0}N$, where $q_0=\varphi(y_0)$, with the following properties:
\begin{enumerate}[(i)]
\item there is a coordinate chart $\Psi:B^N_{r_0}(q_0)\rightarrow T_{q_0}N$ diffeomerphic onto its image, under which $\varphi(B^\Sigma_{r_0}(x_0))$ is the graph of $u$ in the following sense:
\begin{equation*}
\Psi\circ\varphi(B^{\Sigma}_{r_0}(x_0))=\{v+u(v)\vec{N}^\Sigma(y_0):v\in\Omega_0\};
\end{equation*}
\item[(iia)] if $\partial\Sigma\cap B^{\Sigma}_{r_0}(x_0)=\emptyset$, then $y_0=x_0$, ${\varphi}_*B_{\frac{r_0}{2}}(0)\subset\Omega_0$ and $\|u|_{B_{\frac{r_0}{4}}(0)}\|_{2,\alpha}<C_0$;

\item[(iib)] if $\partial\Sigma\cap B^{\Sigma}_{r_0}(x_0)\neq\emptyset$, then $y_0\in\partial\Sigma\cap B^{\Sigma}_{r_0}(x_0)$ with $d^\Sigma(x_0,y_0)=d^\Sigma(x_0,\partial\Sigma)$, $\varphi_*B^+_{\frac{r_0}{2}}(0)\subset\Omega_0\subset \varphi_*T^+_{y_0}\Sigma$ and $\|u|_{B^+_{\frac{r_0}{4}}(0)}\|_{2,\alpha}<C_0$

\end{enumerate}
\end{lemma}
\begin{proof}
Let us first choose the appropriate $r_0>0$ depending only on the geometry of $N$ and the constant $A$. 

Firstly fix $0<\varepsilon_0<\min\{\frac{1}{16},\frac{A}{5}\}$. 
Since $N$ is compact, there exists $r_1(N)>0$ sufficiently small such that for every $q\in\partial N$ and $r<r_1$ there {exist} $C_1(N)>0$ and $\Psi:B^N_{r}(q)\rightarrow T^+_q N $ diffeomorphic onto its image so that defining $g={\left(\Psi^{-1}|_{\Psi(B^N_r(q))}\right)}^*g^N$ and $g_0=g^N_q$ metrics on the open set $\Psi(B^N_{r}(q))$ the following holds:
\begin{enumerate}
\item[(a1)] $\Psi(\partial N\cap B^N_{r}(q))\subset T_q\partial N\subset T^+_q N$;
\item[(b1)] $\|g-g_0\|_{C^k}<\varepsilon_0$, for all $k>0$;
\end{enumerate}
Similarly, there exists $r_2(N)>0$ such that for every $p\in N$ with $d^N(p,\partial N)>\frac{r_1}{2}$ and $r<r_2$ there exists $C_2(N)>0$ and $\Psi:B^N_{r}(p)\rightarrow T_p N $ diffeomorphic onto its image so that defining $g=\Psi_*g^N$ and $g_0=g^N_q$ metrics on the open set $\Psi(B^N_{r}(p))$ the following holds:
\begin{enumerate}
\item[(a2)] $\partial N\cap B^N_{r}(p)=\emptyset$;
\item[(b2)] $\|g-g_0\|_{C^k}<\varepsilon_0$ for all $k>0$;
\end{enumerate}
In each case we take the $C^2$ norm with respect to the constant metric $g_0$.
We make the natural identification $T_v(T_pN)=T_pN$ for any $v\in T_pN$.
Observe that conditions $(b1),(b2)$ hold because the metric is $C^2$.

We pick $r_0<\min\{\frac{r_1}{2},r_2\}$ sufficiently small such that $8Ar_0<\frac{1}{2}$, $(1-16Ar_0)>\frac{1}{2}$ and $(1+32Ar_0)^\frac{1}{2}<\frac{3}{2}$.

Now, let $\varphi:\Sigma\rightarrow N$ be a free boundary isometric immersion and $x_0\in\Sigma$ with $|A^\Sigma(x)|\leq A$ on $B^\Sigma_{r_0}(x_0)$.
Let us focus the proof on the case $\partial \Sigma\cap B^\Sigma_{r_0}(x_0)\neq \emptyset$.

Observe that $x_0$ may not necessarily be a boundary point but $d^\Sigma(x_0,\partial \Sigma)<r_0$ so the nearest point $y_0\in\partial \Sigma$ must be $y_0\in B^\Sigma_{r_0}(x_0)$.
Put $p_0=\varphi(x_0)$ and $q_0=\varphi(y_0)$, we first observe that for any $x,x'\in\Sigma$ we have $d^N(\varphi(x),\varphi(x'))\leq d^\Sigma(x,x')$ thus $\varphi(\bar{B}^\Sigma_{r_0}(x_0))\subset \bar{B}^N_{r_0}(p_0)$ and $d^N(p_0,q_0)<r_0$.
It follows that $B^N_{r_0}(p_0)\subset B^N_{2r_0}(q_0)$ and, since $2r_0<r_1$ there exists $\Psi:B^N_{2r_0}(q_0)\rightarrow T^+_{q_0}N$ satisfying conditions $(a1), (b1)$.
We define $g$ and $g_0$ as above and denote by $\Phi=\Psi\circ\varphi|_{\bar{B}^\Sigma_{r_0}(x_0)}$ isometric immersion of $(\bar{B}^\Sigma_{r_0}(x_0),g^\Sigma)$ into $(T_{q_0}^+N,g)$ with free boundary along $T_{q_0}\partial N$.

Since $\|g-g_0\|_{C^2}<\varepsilon_0$ then $\|g^\Sigma-g^\Sigma_0\|_{C^2}<\varepsilon_0$ with $g^\Sigma_0=\Phi^*g_0$. 
Let us relate the geometric quantities of $\bar{B}^\Sigma_{r_0}(x_0)$ computed with respect to $g^\Sigma$ and $g^\Sigma_0$:
\begin{enumerate}[(I)]
\item $(1-\varepsilon_0)d^\Sigma< d^\Sigma_0 < (1+\varepsilon_0)d^\Sigma$,
\item $\|\vec{N}^{\Sigma}-\vec{N}^{\Sigma}_0\|_{C^1}\leq 3\varepsilon_0$,
\item $|A^\Sigma-A^\Sigma_0|< 5\varepsilon_0$ and $|A^\Sigma_0|<2A$ on $B^\Sigma_{r_0}(x_0)$,
\end{enumerate}
where quantities with subscript $0$ are computed with respect to $g^\Sigma_0$ and $\vec{N}^{\Sigma}_0$ is chosen pointing in the same direction as $\vec{N}^{\Sigma}$.

\begin{claim}
We have that $|g_0(\vec{N}_0^\Sigma(y_0),\vec{N}_0^\Sigma(x))-1|<8Ad^\Sigma_0(x,y_0)$ for all $x\in \bar B^\Sigma_{r_0}(x_0)$.
\end{claim}
Let $\gamma:[0,1]\rightarrow \bar B^\Sigma_{r_0}(x_0)$ be a curve joining $y_0$ to $x$.
We compute
\begin{equation*}
\begin{aligned}
\frac{d}{dt}g_0(\vec{N}_0^\Sigma(y_0),\vec{N}_0^\Sigma(\gamma(t))) & = \nabla^0_{\dot{c}(t)}g_0(\vec{N}_0^\Sigma(y_0),\vec{N}_0^\Sigma(x))\\
																																					& = g_0(\vec{N}_0^\Sigma(y_0),\nabla^0_{\dot{c}(t)}\vec{N}_0^\Sigma(x))\\
																																					& = -g_0(\vec{N}_0^\Sigma(y_0),S^\Sigma_0(\dot{c}(t))),
\end{aligned}
\end{equation*}
where $S_0^\Sigma(v)=-\nabla^0_v \vec{N}_0^\Sigma$ is the shape operator of $\Sigma$ with respect to $g_0$.
It follows that
\begin{equation*}
\left|\frac{d}{dt}g_0(\vec{N}_0^\Sigma(y_0),\vec{N}_0^\Sigma(\gamma(t)))\right|\leq 2A\|\dot{c}(t)\|_{g_0}.
\end{equation*}
Hence,
\begin{equation*}
\begin{aligned}
|g_0(\vec{N}_0^\Sigma(y_0),\vec{N}_0^\Sigma(x))-1| & = \left|\int_0^1\frac{d}{dt}g_0(\vec{N}_0^\Sigma(y_0)),\vec{N}_0^\Sigma(\gamma(t)))dt\right|\\
                                                          & \leq 2A\int_0^1{\|\dot{c}(t)\|_{g_0}dt}= 2A\ell^\Sigma_0(\gamma).
\end{aligned}
\end{equation*}

Since both $y_0,x\in \bar{B}^\Sigma_{r_0}(x_0)$, there exists a path $\gamma$ in $\bar{B}^\Sigma_{r_0}(x_0)$ joining $y_0$ to $x$ such that $\ell^\Sigma(\gamma)\leq 2d^{\Sigma}(y_0,x)$.
Thus $\ell^\Sigma_0(\gamma)\leq (1+\varepsilon_0)\ell^\Sigma(\gamma)\leq(1+\varepsilon_0)2d^{\Sigma}(y_0,x)\leq 2\frac{1+\varepsilon_0}{1-\varepsilon_0}d_0^{\Sigma}(y_0,x)<2\frac{17}{15}d_0^{\Sigma}(y_0,x)<4d_0^{\Sigma}(y_0,x)$.
This concludes the proof of claim $1$.

From the choice of $r_0$, we have that $|g_0(\vec{N}_0^\Sigma(y_0),\vec{N}_0^\Sigma(x))-1|<\frac{1}{2}$, in particular $g_0(\vec{N}_0^\Sigma(y_0),\vec{N}_0^\Sigma(x))>\frac{1}{2}>0$.
It follows that for any $x\in \bar{B}^\Sigma_{r_0}(x_0)$ we have that $\Phi_*T_x\Sigma$ is not perpendicular to $\Phi_*T_{y_0}\Sigma$.
Thus, $\Phi_*T_x\Sigma$ can be written as the graph of a linear map over $\Phi_*T_{y_0}\Sigma$.
A direct application of the inverse function theorem implies that for every $x\in \bar{B}^\Sigma_{r_0}(x_0)$ there exists $\delta_x>0$, an open set $\Omega_x\subset \Phi_*T^+_{y_0}\Sigma$ and a $C^2$ function $u_x:\Omega_x\rightarrow \real$ so that $\Phi(B^\Sigma_{\delta_x}(x))=\{v+u_x(v)\cdot\vec{N}^\Sigma(y_0):v\in \Omega_x\}$.
Since $\bar{B}^\Sigma_{r_0}(x_0)$ is compact, we can find an open domain $\Omega_0\subset \Phi_*T^+_{y_0}\Sigma$ and a $C^2$ function $u:\Omega_0\rightarrow \real$ such that:
\begin{equation*}
\Phi(B^\Sigma_{r_0}(x_0))=\{v+u(v)\cdot\vec{N}_0^\Sigma(y_0):v\in \Omega_0\}.
\end{equation*}
It also follows from the construction of $u$ that $\partial\Sigma\cap B^\Sigma_{r_0}(x_0)$ is the graph of $u$ restricted to $\Omega_0\cap \Phi_* T_{y_0}\partial \Sigma$.
\begin{claim}
We have that $\sup|\nabla u|^2\leq 32Ar_0$.
\end{claim}
At $x=v+u(v)\cdot\vec{N}_0^\Sigma(y_0)$ the normal vector with respect to $g_0$ is given by $\vec{N}_0^\Sigma(x)=\pm\frac{1}{\sqrt{1+|\nabla u|^2}}\left(-\nabla u + \vec{N}_0^\Sigma(y_0)\right)$.
On one hand we have
\begin{equation*}
\|\vec{N}_0^\Sigma(x)-\vec{N}_0^\Sigma(y_0)\|^2_{g_0}=\frac{|\nabla u|^2}{1+|\nabla u|^2}+ (\frac{1}{\sqrt{1+|\nabla u|^2}}\pm 1)^2\geq \frac{|\nabla u|^2}{1+|\nabla u|^2}.
\end{equation*}
On the other hand we have
\begin{equation*}
\|\vec{N}_0^\Sigma(x)-\vec{N}_0^\Sigma(y_0)\|^2_{g_0}=2-2g_0(\vec{N}_0^\Sigma(x),\vec{N}_0^\Sigma(y_0))\leq 16Ar_0
\end{equation*}
Therefore,
\begin{equation*}
|\nabla u|^2\leq \frac{16Ar_0}{1-16Ar_0}<32Ar_0,
\end{equation*}
as long as $1-16Ar_0>\frac{1}{2}$, which proves the claim.

Now, pick $\delta_0>0$ to be the largest radius such that $B^+_{\delta_0}(0)\subset\Omega_0\subset \Phi_*T^+_{y_0}\Sigma$.

\begin{claim}
For any pair $v_1,v_2\in B^+_{\delta_0}(0)$ we have that
$d^\Sigma_{0}(v_1+u(v_1)\vec{N}_0^\Sigma(y_0),v_2+u(v_2)\vec{N}_0^\Sigma(y_0))\leq (1+32Ar_0)^{\frac{1}{2}}\|v_1-v_2\|_{g_0}$ and $\delta_0\geq\frac{r_0}{2}$
\end{claim}
We have $\sigma(t)=v_1 + t(v_2-v_1)\in B_{\delta_0}(0)\subset \Omega_0$ for $t\in[0,1]$, from which follows that
\begin{equation*}
\begin{aligned}
d^\Sigma_{0}(v_1+u(v_1)\vec{N}_0^\Sigma(y_0),v_2+ & u(v_2)  \vec{N}_0^\Sigma(y_0)) \leq \ell^\Sigma_0(\sigma+u(\sigma)\vec{N}_0^\Sigma(y_0)) \\
                                                          & \quad =   \int_0^1(1+|\nabla u(v_1 + t(v_2-v_1))|^2)^{\frac{1}{2}}\|v_1-v_2\|_{g_0}dt
\end{aligned}
\end{equation*} 
and the inequality follows from the previous claim.

For the second part observe that if $\left(1+32Ar_0\right)^{\frac{1}{2}}<\frac{3}{2}$ then it follows $d^\Sigma_{0}({v}+u({v})\vec{N}_0^\Sigma(y_0),0)\leq \frac{3}{2}\|v\|_{g_0}$.
Observe that $B^\Sigma_{r_0}(x_0)$ is not equal to $\Sigma$ because $\partial \Sigma$ is a closed curve but $\partial \Sigma \cap B^\Sigma_{r_0}(x_0)$ is the graph of $u|_{\Omega_0\cap \Phi_* T_{y_0}\partial\Sigma}$.
In particular, there exists $\bar{v}\in B^+_{\delta_0}(0)$ so that $d^{\Sigma}({\bar v}+u({\bar v})\vec{N}_0^\Sigma(y_0),0) > \frac{4}{5}r_0$ and $d^\Sigma_{0}({\bar v}+u({\bar v})\vec{N}_0^\Sigma(y_0),0) > (1-\varepsilon_0)\frac{4}{5}r_0>\frac{3}{4}r_0$.

Suppose that $\delta<\frac{r_0}{2}$, then $\|\bar v\|_{g_0}<\frac{r_0}{2}$ and $d^\Sigma_{0}({\bar v}+u({\bar v})\vec{N}_0^\Sigma(y_0),0)< \frac{3r_0}{4}$ which is a contradiction and it proves the claim.

\begin{claim}
We have that $\nabla u$ restricted to $B^+_{\frac{r_0}{2}}(0)$ is a Lipschitz function of Lipschitz constant less than $\frac{9}{2}A$.
\end{claim}
The second fundamental form is a multiple of the Hessian of $u$, thus $|A_0^\Sigma|\leq 2A$ implies $|Hess u|\leq (1+|\nabla u|^2)^{\frac{1}{2}} 2A < 3A$ whenever $(1+32Ar_0)^{\frac{1}{2}}<\frac{3}{2}$.
It follows that for any $v_1,v_2\in\Omega_0$, we have
\begin{equation*}
|\nabla u(v_1)-\nabla u(v_2)|\leq 3A d^\Sigma_0(v_1+u(v_1)\vec{N}_0^\Sigma(y_0),v_2+u(v_2)\vec{N}_0^\Sigma(y_0)).
\end{equation*}
We get that $d^\Sigma_0(v_1+u(v_1)\vec{N}_0^\Sigma(y_0),v_2+u(v_2)\vec{N}_0^\Sigma(y_0))\leq \frac{3}{2}\|v_1-v_2\|_{g_0}$, thus proving the claim.

We note that there exists a linear transformation $H(v,g,\partial g):T_v\Psi(B^N_{2r_0}(q_0))\rightarrow T_v\Psi(B^N_{2r_0}(q_0))$ depending smoothly on $v$, $g$ and its derivatives $\partial g$ such that for every vector field $X$ on $\Psi(B^N_{2r_0}(q_0))$ we have $\divergent_{g}(X)=\divergent_{g_0}(H(g,\partial g)\cdot X)$.
Furthermore, for each $v\in\Omega_0$ we can write $\vec{N}^\Sigma=H(v,u(v),\nabla u(v))\cdot \vec{N}_0^\Sigma$ such that $H(v,u,\nabla u)$ is a linear map depending smoothly on $u,\nabla u$ and $\|H-Id\|_{C^1}<3\varepsilon_0$.

Finally, put $A(v,u,\nabla u)=\frac{1}{\sqrt{1+|\nabla u|^2}}G(g,\partial g)\cdot H(v,u,\nabla u)$ and note that
\begin{equation*}
\divergent_{g_0}\left(A\cdot\nabla u\right)=\divergent_{g_0}\left(A\cdot \vec{N}_0^\Sigma\right)=\divergent_g\left(\vec{N}^\Sigma\right)=f.
\end{equation*}
We know that $u$ is $C^2$ and $\nabla u$ is uniformly Lipschitz so the coefficients of $A$ are Lipschitz on $\Omega_0$ and by \cite{gilbarg-trudinger}*{Theorem $6.19$}, $u$ is $C^{2,\alpha}$.
In particular $\Phi(B^\Sigma_{r_0}(x_0))$ is $C^{2,\alpha}$ and we can improve $\|\vec{N}_0^\Sigma-\vec{N}^\Sigma\|_{C^{1,\alpha}}<3\varepsilon_0$.

Now, we denote by $\eta=\eta(q_0)\in T_{q_0}^+N$ and $\nu=\nu(y_0)\in\Phi_+T_{y_0}^+\Sigma$ the inward conormal vectors to $\partial N,\partial \Sigma$ at $q_0, y_0$ respectively.
The following holds: $g(\vec{N}^\Sigma,\eta)=0$ along $\Phi(B^\Sigma_{r_0}(x_0))\cap \Phi_*T_{y_0}\partial\Sigma$ and $\nu=\eta$ since $\Phi(B^\Sigma_{r_0}(x_0))$ is free boundary with respect to $g$ and $g,g_0$ coincide at $0$.
Furthermore, along $\Omega_0\cap \Phi_*T_{y_0}\partial\Sigma$ we have
\begin{equation*}
\frac{\partial u}{\partial \nu}=g_0(\nabla u,\nu)=\sqrt{1+|\nabla u|^2}g_0(\vec{N}_0^\Sigma,\eta).
\end{equation*}
From which it follows that $\|\frac{\partial u}{\partial \nu}\|_{C^{1,\alpha}}<(1+32Ar_0)^\frac{1}{2}4\varepsilon_0<\frac{3}{2}\cdot\frac{4}{16}<\frac{1}{2}$.
We may now apply Schauder estimates \cite{agmon-douglis-nirenberg1959}*{Theorem 7.1} to obtain $C_1=C_1(N,A)>0$ such that
\begin{equation*}
\|u|_{B^+_{\frac{r_0}{4}}(0)}\|_{C^{2,\alpha}} < C_1\left(\|u|_{B^+_{\frac{r_0}{2}}(0)}\|_{C^0}+\|f\|_{C^{0,\alpha}}+\|\frac{\partial u}{\partial \nu}\|_{C^{1,\alpha}}\right).
\end{equation*}
We conclude the the proof by observing that all quantities on the right are uniformly bounded by constants depending only on $r_0, A$ and $H_0$ se we can find $C_0=C_0(N,A,H_0)>0$ such that $\|u|_{B^+_{\frac{r_0}{4}}(0)}\|_{C^{2,\alpha}} < C_0$ which finishes the proof.

\end{proof}
\begin{remark}
The case when $\partial\Sigma\cap B^\Sigma_{r_0}(x_0)=\emptyset$ follows the exact same proof.
However, we observe that two situations may occur: $\phi(B^\Sigma_{r_0}(x_0))\cap \partial N=\emptyset$ or $B^\Sigma_{r_0}(x_0)\cap \partial N\neq\emptyset$.
In the first situation we use interior coordinates to set up the proof but in the latter we must use boundary adapted coordinates.
Regardless, in both cases we construct functions defined on an open subset of $\Phi_*T_{x_0}\Sigma$ and the remainder of the proof is similar.
\end{remark}

In the following we make precise the notion of graphical convergence of surfaces (see also \cite{ambrozio-carlotto-sharp2018.3}).
Let $\Sigma\subset N$ be an embedded surface.
Given a vector field $X$ defined on an open set $U\subset N$, $\Omega\subset U\cap\Sigma$ an open set and $\varepsilon>0$ we define $V^X_\varepsilon(\Omega)=\{\Phi^X(q,t)\in N:q\in\Omega,t\in(-\varepsilon,\varepsilon)\}$ where $\Phi^X$ is the flow of $X$.

\begin{definition}\label{definition convergence}
Let $(N,\partial N)$ be a $3$-manifold with boundary, $\Sigma_i$ a sequence of smooth properly embedded surfaces in $N$ with $\partial \Sigma_i\subset \partial N$ and $\Sigma$ a smooth properly embedded surface in $N$ with $\partial\Sigma\subset\partial N$.

We say that $\Sigma_i$ converges locally graphically (relative to $\partial N$) in the $C^{k,\alpha}$ topology to $\Sigma$ 
 at $p$ if for some $r>0$ sufficiently small and some vector field $X$ on $B^N_r(p)$ we have:
\begin{enumerate}[(a)]
\item $X|_{B^N_r(p)\cap \Sigma}=\vec{N}^\Sigma$ and $X|_{B^N_r(p)\cap \partial N}\in T\partial N$;
\item For all $\varepsilon>0$ sufficiently small so that $V^X_\varepsilon(B^N_r(p)\cap \Sigma)$ is an open set, there exists $C^{k,\alpha}$ functions $u_i^1,\ldots u_i^{L_p}$ defined on $B^N_r(p)\cap \Sigma$ when $i$ is sufficiently large such that 
\begin{equation*}
\begin{aligned}
& \Sigma_i\cap V^X_\varepsilon(B^N_r(p)\cap \Sigma)=\Sigma_i^1\cup\ldots\cup\Sigma_i^{L_p}\\
 \text{and}& \\
& \Sigma_i^j=\{\Phi^X(q,u_i^j(q)):q\in B^N_r(p)\cap \Sigma\} \text{ for each } j=1,\ldots,L_p;
\end{aligned}
\end{equation*}
\item $u_i^j$ converges to $0$ in the $C^{k,\alpha}$ topology as $i\rightarrow\infty$ for all $j=1,\ldots,L_p$.
\end{enumerate}
If $L$ is constant for sufficiently large $i$ for all $p\in \Sigma$, we say that the convergence is $L$-sheeted everywhere.

Let us now define convegence when the limit is a properly almost embedded surface.
\begin{definition}
Let $(N,\partial N)$ be a $3$-manifold with boundary, $\Sigma_i$ a sequence of smooth properly embedded surfaces in $N$ with $\partial \Sigma_i\subset \partial N$ and $\Sigma$ a smooth properly almost embedded surface in $N$ with $\partial\Sigma\subset\partial N$.

We say that $\Sigma_i$ converges locally graphically (relative to $\partial N$) in the $C^{k,\alpha}$ topology to $\Sigma$ 
 if the following holds:
\begin{enumerate}[(a)]
\item Whenever $p\in\Sigma\setminus\mathcal{S}$ is an embedded point of $\Sigma$ then $\Sigma\cap\nball_r(p)$ is properly embedded and $\Sigma_i\cap\nball_r(p)$ converges to $\Sigma\cap\nball_r(p)$ as in Definition \ref{definition convergence};

\item Whenever $p\in\mathcal{S}$ is a self-touching point of $\Sigma$ then $\Sigma\cap\nball_r(p)=W_1\cup\ldots\cup W_{l_p}$ (as in Definition \ref{almost embedded}) and, for $i$ sufficiently large, $\Sigma_i\cap\nball_r(p)=\Sigma_{i,1}\cup\ldots\cup\Sigma_{i,l_p}$ so that $\Sigma_{i,m}$ converges to the embedded component $W_m$ as in Definition \ref{definition convergence} for all $m=1,\ldots,l_p$.
\end{enumerate}
\end{definition}

\end{definition}

\begin{remark}
Firstly observe that if $\Sigma_i$ converges to $\Sigma$ on $B^N_r(p)$ and $V^X_\varepsilon$ relative to some extension $X$ of $\vec{N}^\Sigma$ then it also converges on $B^N_{r'}(p)$ and $V^{X'}_{\varepsilon'}$ for any other extension $X'$ of $\vec{N}^\Sigma$ and some $r'<r$, $\varepsilon'>0$.
Secondly, let $U\subset N$ be an open set with $U\cap \Sigma$ is connected and $\Sigma_i$ converges graphically to $\Sigma$ at all points of $U\cap \Sigma$, then we can find $X$ an extension of $\vec{N}^\Sigma$ on $U$ such that $\Sigma_i\cap V^X_\varepsilon(U\cap\Sigma)$ is given as the graph of functions $u_i^j$ defined on $U\cap \Sigma$, $L_i=L$ is constant for $i$ sufficiently large and $u_i^j$ converges to $0$ in the $C^{k,\alpha}$ topology on compact sets of $U\cap\Sigma$.
\end{remark}


\section{Curvature estimate}

In this section we prove an improvement from uniformly small total curvature estimate to uniform local pointwise curvature estimate.
The proof is inspired by \cite{white1987}*{Theorem 1} and we focus on the local estimates around a boundary point.
We point out as well that the same proof holds even if the surface is not CMC but has uniformly bounded $C^{0,\alpha}$ mean curvature.
We also prove curvature estimates for weakly stable free boundary CMC surfaces.

\begin{theorem}\label{curvature estimates}
Let $N$ be a compact Riemannian $3$-manifold with smooth boundary. There exists a small enough $r_0>0$ such that the following holds: whenever $\Sigma$ is a properly immersed CMC surface in $N$, $Q\in\Sigma$, $\partial\Sigma\cap \nball_{r_0}(Q)$ is either empty or free boundary in $\partial N\cap \nball_{r_0}(Q)$ and the mean curvature of $\Sigma$ satisfies $H_{\Sigma}\leq H_0$,
then there exist $\varepsilon_0>0$ depending on $\nball_{r_0}(Q)$ and $H_0$ such that if $\int_{\Sigma\cap \nball_{r_0}(Q)} |A^{\Sigma}|^2\leq \varepsilon_0$, then
\begin{equation*}
\max_{0\leq \sigma\leq r_0}\left(\sigma^2\sup_{\Sigma\cap \nball_{r_0-\sigma}(Q)}|A^\Sigma|^2\right)\leq C_0
\end{equation*}
where the constant $C_0$ only depends on geometry of $\nball_{r_0}(Q)$ and $H_0$.
\end{theorem}
\begin{proof}
Suppose false, that is, for $r_n\rightarrow 0$ and $\varepsilon_n\rightarrow 0$ there exist free boundary CMC surfaces $\Sigma_n\subset N$ and $Q_n\in\Sigma_n$ satisfying:
\begin{enumerate}[(i)]
\item $H_n\leq H_0$;
\item $\int_{\Sigma_n\cap \nball_{r_n}(Q_n)} |A_n|^2\leq \varepsilon_n$ and
\item $\max_{0\leq \sigma\leq r_n}\left(\sigma^2\sup_{\Sigma\cap \nball_{r_n-\sigma}(Q_n)}|A_n|^2\right) > n$.
\end{enumerate}
Pick $0<\sigma_n< r_n$ such that
\begin{equation*}
\sigma_n^2\sup_{\Sigma\cap \nball_{r_n-\sigma_n}(Q_n)}|A_n|^2=\max_{0\leq \sigma\leq r_n}\left(\sigma^2\sup_{\Sigma\cap \nball_{r_n-\sigma}(Q_n)}|A_n|^2\right)
\end{equation*}
and write $\lambda_n^2 = \sup_{\Sigma\cap \nball_{r_n-\sigma_n}(Q_n)}|A_n|^2$.
For each $n$ there exists $z_n\in\Sigma\cap \nball_{r_n-\sigma_n}(Q_n)$ such that $|A_n(z_n)|>\frac{\lambda_n}{2}$.

By taking a subsequence we have $Q_n\rightarrow Q$ and $\nball_{r_n}(Q_n)$ contained in a geodesic ball of $N$.
Without loss of generality we may assume that either $\Sigma_n'=\Sigma_n\cap \nball_{r_n}(Q_n)\subset B_{r_n}(0)\subset \real^3$ or $\Sigma_n'\subset B^+_{r_n}(0)\subset \real^2\times \real_+$ in case $Q_n$ is a boundary point for all $n$ sufficiently large, in both cases we induce a metric $g$ on the corresponding subset of $\real^3$.
Since the following arguments are exactly the same in both cases, we will omit the indication of the half-ball.

We denote by $F_n:B_{r_n}(0)\rightarrow \real^3$ the map $F_n(q)=\lambda_n(q-z_n)$ and we define $\tilde\Sigma_n=F_n(\Sigma_n')$, which by abuse of notation we write as $\lambda_n(\Sigma_n'-z_n)$.
Note that $\tilde\Sigma_n$ is a properly immersed surface in $F_n(B_{r_n}(0))$ with the metric $g_n={F_n}_*g$.
We further observe that $B_1(0)\subset F_n(B_{r_n}(0))$ for $n$ sufficiently large and $g_n$ converges to the Euclidean metric.


Furthermore, $\tilde\Sigma_n$ satisfies:
\begin{enumerate}[(a)]
\item $|\tilde A_n(\tilde x_n)|\leq 2$ for all $\tilde x_n\in\tilde \Sigma_n\cap B_1(0)$ and $n$ sufficiently large;
\item $|\tilde A_n(0)|>\frac{1}{2}$ and
\item $\int_{\tilde\Sigma_n} |\tilde A_n|^2\leq \varepsilon_n$
\end{enumerate}
Indeed, if $\tilde x_n\in B_1(0)$ then $\tilde x_n = \lambda_n(x_n-z_n)$ with $x_n\in \Sigma_n\cap B_{\frac{1}{\lambda_n}}(z_n)$.
It follows that $x_n\in \Sigma_n\cap B_{r_n-(\sigma_n-\frac{1}{\lambda_n})}(Q_n)$.
Since 
\begin{equation*}
\left(\sigma_n-\frac{1}{\lambda_n}\right)^2\sup_{\Sigma\cap \nball_{r_n-(\sigma_n-\frac{1}{\lambda_n})}(Q_n)}|A_n|^2\leq \max_{0\leq \sigma\leq r_n}\left(\sigma^2\sup_{\Sigma\cap \nball_{r_n-\sigma}(Q_n)}|A_n|^2\right)=\sigma_n^2\lambda_n^2,
\end{equation*}
it implies that 
\begin{equation*}
|A_n(x_n)|^2<\left(\frac{1}{1-\frac{1}{\sigma_n\lambda_n}}\right)^2\lambda_n^2.
\end{equation*} 
We know that $\sigma_n^2\lambda_n^2 > n$ thus $|A_n(x_n)|<2\lambda_n$ for $n$ sufficiently large, which proves (a).
Property (b) follows from rescaling and (c) holds because the total curvature is scale invariant.

Let $\hat \Sigma_n$ denote the connected component of $\tilde \Sigma_n\cap B_1(0)$ containing $0$.
We may further assume, after an ambient rotation and translation that $T_{0}\hat \Sigma_n=\{x_3=0\}$.

Using property (a), it follows from Lemma \ref{schauder estimates}(i),(ii) that there exists $\hat r_0$ independent of $n$ such that, under the appropriate identifications, $\hat \Sigma_n\cap B_{\hat r_0}(0)$ is the graph of a function $u_n$ satisfying $\|u_n|_{B_{\frac{\hat r_0}{4}}(0)}\|_{2,\alpha}<\hat C_0$, where $\hat C_0$ is independent of $n$.
If $0$ is a boundary point then the domain of $u_n$ contains a half-ball instead but Lemma \ref{schauder estimates}(i),(iii) apply similarly.

Finally, $u_n$ converges, up to a subsequence, in the $C^2$ topology to a function $u_\infty$.
If we denote its graph by $\hat{\Sigma}_\infty$ then it satisfies $|\hat A_\infty(0)|\geq\frac{1}{2}$ from (b) and $\int_{\hat\Sigma_\infty}|\hat A_\infty|^2=0$ from (c), which is a contradiction and completes the proof of the Theorem.
\end{proof}

\begin{definition}
Let $\Sigma$ be a free boundary CMC surface in $N$ and $U\subset N$ be an open set.
We say that $\Sigma$ is weakly stable in $U$ if for every compactly supported function $\phi$ in $\Sigma\cap U$ with $\int_{\Sigma\cap U}\phi d\volume_{\Sigma}=0$ we have
\begin{equation*}
\frac{d^2}{dt^2}|_{t=0}\area(\Sigma(t))\geq 0,
\end{equation*}
where $\Sigma(t)$ is a variation of $\Sigma$ with respect to $\phi$.
\end{definition}

To prove the multiplicity $1$ convergence we will need curvature estimates for free boundary CMC surfaces.
In the case of surfaces the proof is simple and it follows directly from a Bernstein-type theorem \cite{hong-saturnino2023}*{Theorem $1.5$} and a blow-up argument as above.
We include the proof here for completeness.

\begin{theorem}\label{curvature estimates 2}
Let $N$ be a compact Riemannian $3$-manifold with smooth boundary, $Q\in N$ and $H_0>0$.
There exista $r_0>0$ and $C_0=C_0(H_0,\nball_{r_0}(Q))>0$ such that the following holds: whenever $\Sigma$ is a smooth properly immersed free boundary CMC surface in $N$, $Q\in\Sigma$, $H_{\Sigma}\leq H_0$ and $\Sigma$ is weakly stable in $N$, then
\begin{equation*}
\max_{0\leq \sigma\leq r_0}\left(\sigma^2\sup_{\Sigma\cap \nball_{r_0-\sigma}(Q)}|A^\Sigma|^2\right)\leq C_0
\end{equation*}
\end{theorem}
\begin{proof}
Suppose false, that is, for $r_n\rightarrow 0$ there exist free boundary CMC surfaces $\Sigma_n\subset N$ satisfying:
\begin{enumerate}[(i)]
\item $H_{\Sigma_n}\leq H_0$;
\item $\max_{0\leq \sigma\leq r_n}\left(\sigma^2\sup_{\Sigma_n\cap \nball_{r_n-\sigma}(Q)}|A_n|^2\right) > n$.
\item $\Sigma_n$ is weakly stable in $N$
\end{enumerate}
Pick $0<\sigma_n< r_n$ such that
\begin{equation*}
\sigma_n^2\sup_{\Sigma\cap \nball_{r_n-\sigma_n}(Q)}|A_n|^2=\max_{0\leq \sigma\leq r_n}\left(\sigma^2\sup_{\Sigma_n\cap \nball_{r_n-\sigma}(Q_n)}|A_n|^2\right)
\end{equation*}
and write $\lambda_n^2 = \sup_{\Sigma\cap \nball_{r_n-\sigma_n}(Q)}|A_n|^2$.
For each $n$ there exists $z_n\in\Sigma\cap \nball_{r_n-\sigma_n}(Q)$ such that $|A_n(z_n)|>\frac{\lambda_n}{2}$.

If $Q\in N\setminus\partial N$ we may assume that $r_n<\min\{\frac{d(Q,\partial N)}{2},\textup{inj}(Q)\}$, where $\textup{inj}(Q)$ is the injectivity radius of $N$ at $Q$.
If $Q\in\partial N$ we may assume that $\nball_{r_n}(Q)$ is a geodesic half-ball adapted to the boundary of $N$.

We assume $Q\in\partial N$ since the interior case is similar.
For $n$ sufficently large we may further assume that $\nball_{r_n}(Q)$ is diffeomorphic to $B^+_{r_n}(0)\subset T^+_Q N$.
If we idenfity $T^+_Q N=\real^2\times \real_+$ we may further assume that $\partial N\cap \nball_{r_n}(Q)$ is mapped onto $\{x\in \real^2\times \real_+: x_3=0\}$
We will denote by $\Sigma_n'$ the identification of $\Sigma_n\cap\nball_{r_n}(Q)$ in $B^+_{r_n}(0)$ and by abuse of notation the induced metric $g$ on $B^+_{r_n}(0)$.

Let $F_n:B^+_{r_n}(0)\rightarrow \real^3$ be the map $F_n(q)=\lambda_n(q-z_n)$ and define $\tilde\Sigma_n=F_n(\Sigma_n')$.
Note that $\tilde\Sigma_n$ is a properly immersed free boundary CMC surface in $F_n(B^+_{r_n}(0))$ with the metric $g_n={F_n}_*g$.

Observe that $\tilde\Sigma_n$ satisfies:
\begin{enumerate}[(a)]
\item If $R<\frac{\sqrt{n}}{2}$ then $|\tilde A_n(\tilde x)|\leq 4$ for all $\tilde x\in\tilde \Sigma_n\cap B^+_R(0)$;
\item $|\tilde A_n(0)|>\frac{1}{2}$;
\item $H_{\tilde{\Sigma_n}}\leq\lambda_n^{-1}H_0$ and
\item $\tilde\Sigma_n$ is weakly stable in $F_n(B^+_{r_n}(0))$.
\end{enumerate}
Indeed, if $\tilde x\in B^+_R(0)$ then $\tilde x = \lambda_n(x-z_n)$ with $x\in \Sigma_n'\cap B_{\frac{R}{\lambda_n}}(z_n)$.
Since $\sigma_n\lambda_n > 2R$ it follows that $x\in \Sigma_n'\cap B_{r_n-\frac{\sigma_n}{2}}(Q)$ and
\begin{equation*}
\left(\frac{\sigma_n}{2}\right)^2\sup_{\Sigma\cap \nball_{r_n-\frac{\sigma_n}{2}}(Q)}|A_n|^2\leq \max_{0\leq \sigma\leq r_n}\left(\sigma^2\sup_{\Sigma\cap \nball_{r_n-\sigma}(Q_n)}|A_n|^2\right)=\sigma_n^2\lambda_n^2,
\end{equation*}
it implies that 
\begin{equation*}
|A_n(x_n)|^2<4\lambda_n^2,
\end{equation*} 
which proves (a).
Properties (b) and (c) follow from rescaling.
To prove property (d) let $\phi$ be a compactly supported function on $\tilde\Sigma_n$ satisfying $\int_{\tilde\Sigma_n}\phi d\volume_{g_n}=0$.
Then $\phi_n(q)=\phi(F_n(q))$ is a compactly supported function on $\Sigma_n'$ satisfying $\int_{\Sigma_n'}\phi_n d\volume_g=\lambda_n^{-2}\int_{\tilde\Sigma_n}\phi d\volume_{g_n}=0$.
If $\tilde\Sigma_n(t)$ is a variation of $\tilde\Sigma_n$ with respect to $\phi$, then $\Sigma_n'(t)=F_n^{-1}(\tilde\Sigma_n(t))$ is a variation of $\Sigma_n'$ with respect to $\phi_n$ for $t$ sufficiently small.
We observe that $\frac{d^2}{dt^2}|_{t=0}\area(\tilde\Sigma_n(t))=\lambda_n ^{-2}\frac{d^2}{dt^2}|_{t=0}\area(\Sigma_n'(t))>0$.
That is, $\tilde\Sigma_n$ is weakly stable in $F_n(B^+_{r_n}(0))$.

Let $\hat \Sigma_n$ denote the connected component of $\tilde \Sigma_n$ containing $0$.
Observe that $F_n(B^+_{r_n}(0))$ converge to $\real^2\times\real_+$ and $g_n$ converge to the Euclidean metric smoothly on compact sets.
Fix $R>0$, and let $\delta_0=r_0(N,4)$ as given by Lemma \ref{schauder estimates} with a fixed $\alpha\in(0,1)$.
It follows from property $(a)$ that for all $n>4R^2$ and $\hat x\in\nball_R$, $\hat\Sigma_n\cap B_{\delta_0}(\hat x)$ can be written as a graph of a function with uniform $C^{2,\alpha}$ estimates depending only on $N$, $H_0$ and $\alpha$.
By covering $B^+_R$ with finitely many balls $B_{\delta_0}$ where the number of balls depend only on $R$ and applying a diagonal argument we can construct a subsequence of the immersions of $\hat\Sigma_n$ in $F_n(B^+_{r_n}(0))$ that converge to an immersion in $\real^2\times\real_+$ in the $C^{2,\beta}$ topology for all $0<\beta<\alpha$.
If we take $R_n>0$ tending to infinity and a diagonal sequence from the argument above, we obtain a $C^{2,\beta}$ surface $\Sigma_\infty$ in $\real^2\times\real_+$ and subsequence of $\hat\Sigma_n$ that converges to $\Sigma_\infty$ in the $C^{2,\beta}$ topology on compact sets (see also \cite{smith.g2007}*{Theorem 1.2}).

Properties $(c)$, $(d)$ and $C^2$ convergence imply that $\Sigma_\infty$ is a weakly stable free boundary minimal surface in $\real^2\times\real_+$ with respect to the Euclidean metric and boundary on $\{x\in\real^2\times\real_+:x_3=0\}$.
It follows from \cite{hong-saturnino2023}*{Theorem 1.5} that $\Sigma_\infty$ is a half plane.
However, property $(b)$ implies that $|A_{\Sigma_\infty}|(0)\geq\frac{1}{2}$ which is a contradiction and concludes the proof of the Theorem.

\end{proof}


\section{Removable singularities}

As we will see later, the compactness result does not give us smooth convergence everywhere.
The points in which we do not have sufficient curvature estimates are potential singularities either because of a neckpinching phenomenum where the curvature may blow up, or self-touching points where the convergence is not single sheeted.
However, we are still able to prove that these are removable singularities so the limiting surface is still a smooth object.

We are going to prove that if the total curvature is bounded on a CMC surface, then isolated singularities are removable.
This is an adaptation of \cite{white1987}*{Theorem 2} and the arguments are the same except for the foliation argument to prove uniqueness of the tangent cone.
We are going to focus on the case in which the singularity is along the boundary, but the same result holds for interior singularities and the proof follows the exact same arguments.

The idea of the proof is to improve the integral curvature bound to a pointwise curvature decay near the singularity to show that the tangent cones are totally geodesic.
By adapting the foliation argument of White \cite{white1987} we prove that the tangent cone is unique from which we can show that near the singularity the surface is indeed the graph of a $C^1$ function.
We can then improve it further using elliptic regularity.

Firstly, let us prove the existence of a local CMC foliation with free boundary which will be needed later.
This is a straightforward adaptation of \cite{ambrozio-carlotto-sharp2018.3}*{Section $3$} together with White's approach to deal with a family of functionals \cite{white1987}*{Appendix}.

Let $\theta\in(0,\frac{\pi}{4})$ and define $D_\theta = \{ x \in\real^2 : (x_1+a)^2+x_2^2 \leq 1 \text{ and } x_1\geq 0\}$, where $a=\cos^{-1}(\theta)\in(\frac{1}{\sqrt{2}},1)$.
This is the part of the disk of radius $1$ centered on the $x_1$-axis that intersects the line $x_1=0$ at angle $\theta$.
Its boundary components are denoted by $\partial_0D_\theta=\partial D_\theta\cap\{x\in\real^2:x_1=0\}$ and $\partial_+D_\theta=\overline{ \partial{D_\theta\setminus \partial_0D_\theta}}$.
The regular cylinder over $D_\theta$ in $\real^3$ is denoted by $C_\theta=D_\theta\times\real$, with corresponding boundary components $\partial_0C_\theta=\partial_0 D_\theta\times\real$ and $\partial_+C_\theta=\partial_+ D_\theta\times\real$.


\begin{figure}[H]
\includegraphics[scale=0.5]{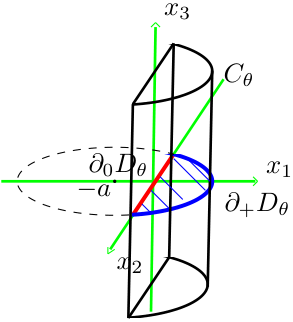}
\end{figure}

Given a function $f\in C^{2,\alpha}(D_\theta)$, we define $N_g^+(f)$ to be the normal vector over $\graph(f)$ with respect to $g$ pointing in the positive direction of the $x_3$-axis, that is, $g(N_g^+(f),\frac{\partial}{\partial x_3})>0$.
We write $H_g^+(f)=g(\vec{H}_g(f),N_g^+(f))$ as the scalar mean curvature with respect to $N_g^+(f)$.
In particular, $\vec{H}_g(f)$ points in the positive direction of the $x_3$-axis when $H_g^+(f) > 0$ and in the negative direction otherwise.

Let us denote by $X$ the space of $C^{2,\alpha}$ metrics on $C_\theta$ and define the map

\begin{equation*}
\Phi:\real\times\real\times X \times C^{2,\alpha}(\partial_+D_\theta) \times C^{2,\alpha}(D_\theta)\rightarrow C^{0,\alpha}(D_\theta)\times C^{1,\alpha}(\partial_0 D_\theta)\times C^{2,\alpha}(\partial_+ D_\theta)
\end{equation*}
by
\begin{equation*}
\Phi(h,t,g,w,u)=\left(H^+_g(t+u)-h,\frac{\partial}{\partial\eta_g}(t+u),u_{|_{\partial_{+}D_{\theta}}}-w\right),
\end{equation*}
where $\eta_g$ is the inward conormal vector along $\partial_0D_\theta$.

\begin{proposition}[{\cite{ambrozio-carlotto-sharp2018.3}*{Proposition $21$}}]\label{proposition foliation}
For every $t_0\in\mathbb{R}$, there exist a neighbourhood $U_{t_0}$ of the Euclidean metric $\delta$ in $X$, $\varepsilon_{t_0}>0$ and
\begin{equation*}
u:(-\varepsilon_{t_0},\varepsilon_{t_0})\times(t_0-\varepsilon_{t_0},t_0+\varepsilon_{t_0})\times U_{t_0}\times B^{C^{2,\alpha}(\partial_+D_\theta)}_{\varepsilon_{t_0}}(0)\rightarrow C^{2,\alpha}(D_\theta)
\end{equation*}
so that $t\mapsto \graph(t+u(h,t,g,w))$ defines a $C^{2,\alpha}$ foliation of $D_\theta\times[t_0-\frac{\varepsilon_{t_0}}{2},t_0+\frac{\varepsilon_{t_0}}{2}]$ by surfaces with constant mean curvature $h$ with respect to the metric $g$, free boundary along $\partial_0 C_\theta$ and $(t+u)_{|_{\partial_+D_\theta}}=t+w$.

Furthermore, if $h>0$ then $\vec{H}_g(t+u)$ points in the positive direction of the $x_3$-axis and in the negative direction when $h<0$.
 \end{proposition}
 
\begin{proof}
Observe that $\Phi$ defined above is a $C^1$ function and $D_5\Phi(0,t_0,\delta,0,0)$ defines the same isomorphism as in \cite{ambrozio-carlotto-sharp2018.3}*{Appendix B}.
The result then follows from the Implicit Function Theorem.
\end{proof}

Let us denote by $B^+_1=\{x\in\real^3 : \|x\|\leq 1, x_1\geq 0\}$ the upper half-ball in $\real^3$ and $\partial_0B^+_1=\partial B^+_1\cap\{x\in\real^3:x_1=0\}$.

\begin{theorem}\label{removable singularity boundary}
Let $g$ be a Riemannian metric on $B^+_1$ and $\Sigma$ be a smooth, properly embedded, CMC surface in $B^+_1\setminus\{0\}$, $\partial \Sigma\subset \partial B^+_1$, free boundary on $\partial_0 B^+_1\setminus\{0\}$ and $0\in\overline{\partial\Sigma}$.
Suppose $\int_\Sigma|A^\Sigma|^2\leq C$ then $\Sigma\cup\{0\}$ is a smooth, properly almost embedded, CMC surface in $B_1^+$.
\end{theorem}

\begin{proof}
Let $r_0>0$ and $\varepsilon_0>0$ be as in Theorem \ref{curvature estimates}.
Pick $\delta>0$ sufficiently small so that $\int_{\Sigma\cap B^+_\delta}|A^\Sigma|^2\leq \varepsilon_0$ where $B_\delta^+=\delta B_1^+$.
It follows from Theorem \ref{curvature estimates} that
\begin{equation*}
|A^\Sigma(x)|d_g(x,0)\leq C_0,
\end{equation*}
whenever $d_g(x,0)<\delta$.

\begin{claim}
Every tangent cone of $\Sigma$ at $0$ is a union of half-planes in $\real^3\setminus\{0\}$.
\end{claim}
Let $r_i\rightarrow\infty$ be any sequence and define $F_i:B_1^+\rightarrow \real ^3$ as $F_i(x)=r_ix$, $\Sigma_i=F_i(\Sigma)$, which by abuse of notation we write $\Sigma_i=r_i\Sigma$.
It follows that $\Sigma_i$ is a properly embedded CMC surface in $(F_i(B_1^+\setminus\{0\}),({F_i^{-1}|_{F_i(B_1^+\setminus\{0\})})^*g})$, $F_i(B_1^+)$ converges to $T_0^+B_1^+$ which we identify with $\real^3_+=\{x\in\real^3:x_1\geq 0\}$, which implies that for any compact set $K\subset \real^3_+$ we have $K\subset F_i(B_1^+)$ for $i$ sufficiently large.
Furthermore, ${({F_i^{-1}|_{F_i(B_1^+\setminus\{0\})}})^*g}$ converges to the Euclidean metric on compact sets.

Observe that the curvature estimate above is scale invariant, so $\Sigma_i$ satisfies the same curvature bounds whenever $d_g(x,0)<r_i\delta$.
It follows that, up to a subsequence, $\Sigma_i$ converges locally graphically in the $C^{1,\alpha}$ topology on compact sets to a complete surface $\Sigma_\infty$ in $T^+_0B_1^+=\real^3_+$ equiped with the Euclidean metric.
Lemma \ref{schauder estimates} implies that $\Sigma_i$ in fact converges locally graphically in the $C^{2,\alpha}$ topology on compact sets of $\real^3_+\setminus\{0\}$.
In particular, $\Sigma_\infty$ has {non-empty} free boundary on $\{x\in\real^3:x_1=0\}\setminus\{0\}$ and for any compact set $K\subset\real^3_+\setminus\{0\}$ we have
\begin{equation*}
\int_{K\cap\Sigma_\infty}|A_{\infty}|^2 = \lim_{i\rightarrow\infty}\int_{K\cap \Sigma_i}|A_i|^2= \lim_{i\rightarrow\infty}\int_{(r_i^{-1}K)\cap \Sigma}|A^\Sigma|^2=0.
\end{equation*}
That is, $A_\infty=0$. 
Hence $\Sigma_\infty$ is an union of half-planes perpendicular to $\{x\in\real^3:x_1=0\}$.

\begin{claim}
If $\delta>0$ is sufficently small then $\Sigma\cap B^+_\delta\setminus\{0\}$ is topologically a finite and disjoint union of disks, half-disks or half-disks punctured at $0$ with free boundary on $\partial_0 B^+_\delta\setminus\{0\}$.
\end{claim}
Firstly, we improve the curvature estimates.
Fix $y\in\Sigma_i\cap B^+_{r_i\delta}$ and put $x=r_i^{-1}y\in\Sigma\cap B^+_{\delta}$.
Then $|A^\Sigma(x)|d_g(x,0)=|A_i(y)|d_g(y,0)\rightarrow 0$ as $i\rightarrow\infty$ from the previous claim.
Thus $\lim_{x\rightarrow 0}|A^\Sigma(x)|d_g(x,0)=0$.

Secondly, we use a standard Morse Theory argument.
Let $(\tilde U, \tilde g)$ be a Riemannian prolongation of $(B_\delta^+,g)$ such that $\tilde U \cap B^+_1=B^+_\delta$, $\tilde g= g$ on $\bar{B}_\delta^+$ and $\partial_0 B_1^+\cap \tilde U$ separates $\tilde U$ into two connected open sets one of which is $B_\delta^+$.
If $\delta >0$ is sufficiently small we may further assume that $\partial \tilde U$ is convex.

Let $\tilde f:\tilde U\rightarrow \real $ be defined as $\tilde f(x)=\frac{1}{2}d_{\tilde g}^{\tilde U}(x,0)^2$.
If $x\in \Sigma$ and $v\in T_x\Sigma$ then, 
{
\begin{equation*}
\textup{Hess}\tilde{f}(v,v)_x=\tilde{g}(\vec{A}^\Sigma(v,v)_x,\vec{N}^\Sigma(x) )_x\tilde{g}(\vec{N}^\Sigma(x),\dot\gamma_{\tilde g}(1))_x+\tilde{Q}_x(v,v)
\end{equation*}
where $\vec{A}^\Sigma(v,w)_x=\nabla^{\tilde{g}}_vw-\nabla^\Sigma_vw$,
}
$\tilde{Q}_x(v,v)\geq\frac{1}{2}\|v\|^2_{\tilde{g}}$ for $\delta$ sufficiently small and  $\gamma_{\tilde g}:[0,1]\rightarrow \tilde U$ is the unique minimizing geodesic in $\tilde U$ with respect to $\tilde g$ with $\gamma_{\tilde g}(0)=0$ and $\gamma_{\tilde g}(1)=x$.
This follows from Gauss Lemma and writing $\nabla^{\tilde g} f = \nabla^\Sigma f + \tilde{g}(\vec{N}^\Sigma, \nabla^{\tilde g} f)\vec{N}^\Sigma$.
Observe that for any $x\in B_\delta^+$ we have $d_{\tilde g}^{\tilde U}(x,0)= d_g(x,0)$ and $\|\dot\gamma_{\tilde g}(t)\|_{\tilde{g}}=d_{\tilde g}^{\tilde U}(\gamma_{\tilde g}({1}),0)$ so
\begin{equation*}
\begin{aligned}
|{\tilde{g}(\vec{A}^\Sigma(v,v)_x,\vec{N}^\Sigma(x) )_x}\tilde{g}(\vec{N}^\Sigma(x),\dot\gamma_{\tilde g}(1))_x| & \leq|A^\Sigma(x)|d_{\tilde g}^{\tilde U}(x,0)\|v\|_{\tilde g}\\
                                                                                                                             & \leq |A^\Sigma(x)|d_{g}(x,0)\|v\|_{\tilde g}.
\end{aligned}
\end{equation*}

Therefore {$\tilde{g}(\vec{A}^\Sigma(v,v)_x,\vec{N}^\Sigma(x) )_x \tilde{g}(\vec{N}^\Sigma(x),\dot\gamma_{\tilde g}(1))_x$} tends to $0$ as $x$ tends to $0$.
In particular, for $\delta$ sufficiently small any critical point of $f=\tilde{f}|_\Sigma$ in a small neighbourhood of $0$ is a strict local minimum, even if the critical point is along the boundary. 
It follows from Morse Theory for manifolds with boundary \cite{jankowski-rubinsztein1972} that every connected component of $\Sigma\cap B^+_\delta$ is a disk with a single critical point, a free boundary half-disk with a single critical point or a free boundary half-disk punctured at $0$ without critical points.
Since critical points cannot accumulate, it must be finite for $\delta$ sufficiently small.
This proves the claim.

Next, we shall prove that for each punctured half-disk its tangent cone is unique.
Pick $\hat\Sigma\subset\Sigma\cap B^+_\delta\setminus\{0\}$ a connected component that corresponds to one of the punctured half-disks.
Let $\hat\Sigma_i$ be the corresponding dilation by $r_i\rightarrow\infty$ and $\hat\Sigma_\infty\subset\Sigma_\infty$ its tangent cone at $0$.
As we have seen above, $\hat\Sigma_\infty$ is a half-plane perpendicular to $\{x\in\real^3:x_1=0\}$.

{
Without loss of generality let us identify $\hat\Sigma_\infty$ with the half-plane $P_+=\{x\in\real^3: x_3=0\}\cap\{x\in\real^3: x_1\geq 0\}$.
We may also assume, possibly taking another subsequence, that the mean curvature vector of $\hat\Sigma_i$ points in the positive direction of the $x_3$-axis for $i$ sufficiently large. 
Let $D_\theta$ and $C_\theta$ be as in Proposition \ref{proposition foliation}. 
\begin{claim}
For any $R>0$ there exists $i_0$ sufficiently large such that $\hat\Sigma_i\cap\partial_+D_\theta\times(-R,R)$ is the graph of a single function $w_i=w_{r_i}$ over $\partial_+D_{\theta}$ for all $i>i_0$.
Furthermore, $\|w_i\|_{2,\alpha}\rightarrow 0$.
\end{claim}

Since $\hat\Sigma_i\cap D_\theta\times(-R,R)$ converges grahically to $P_+\cap D_\theta\times(-R,R)$ away from $0$ and $\hat{\Sigma}_i$ is obtained by dilation, we may take $i_0$ sufficiently large so that every connected component of $\hat\Sigma_i\cap D_\theta\times(-R,R)$ contains $0$ in its closure for all $i>i_0$.
Indeed the distance between $0$ and each connected component of $\hat\Sigma_i\cap C_\theta$ that does not contain $0$ in its closure can only increase by dilation, so for $i$ sufficiently large we may assume that it is larger than $4R$ which is sufficently to conclude that such component does not belong to $\cap D_\theta\times(-R,R)$.

Now, suppose there were two connected components of $\hat\Sigma_i\cap D_\theta\times(-R,R)\setminus\{0\}$ containing $0$ in its closure for some $i>i_0$ (recall $0\in\bar\Sigma\setminus \Sigma$).
Two such components, namely $E_i^1,E_i^2$ must correspond to embedded punctured half-disks and the corresponding disjoint sets $\frac{1}{r_i}E_i^1,\frac{1}{r_i}E_i^2$ belong to $\hat\Sigma$.
This is a contradiction since $\hat\Sigma$ is a connected component corresponding to a half-disk with a single puncture in which the distance function $f$ defined above has no critical points.
Indeed the existence of two such components would imply that $\hat\Sigma$ corresponds topologically to a half-disk with at least two punctures.
Equivalently, it would imply the existence of a critical point of $f|_{\hat\Sigma}$.

The final conclusion of the claim follows because $\hat\Sigma_i$ converges graphically to the half-plane $P_+$ away from $0$.
}

\begin{claim}
The tangent cone $\hat\Sigma_\infty$ is unique, that is, it is independent of the blow-up sequence $r_i\rightarrow\infty$.
\end{claim}

{Let $t_0=0$ and and $\varepsilon_0>0$ as in Proposition \ref{proposition foliation}.
In view of the previous claim, for $R_0=2\varepsilon_0$ we know that $\hat\Sigma_i\cap\partial_+D_{\theta}\times (-R_0,R_0)$ is the graph of a function $w_i=w_{r_i}$ over $\partial_+D_{\theta}$ for $i$ sufficiently large and $\|w_i\|_{2,\alpha}\rightarrow 0$.
We further assume that $i$ is large enough so that $\|w_i\|_{2,\alpha}<\varepsilon_0$.
}

The mean curvature of $\hat\Sigma_i$ is given by $\hat{H}_i={r_i}^{-1}H_\Sigma$ which is constant at each $i$ and tends to $0$ as $i$ tends to infinity.
If we denote by $g_i$ the corresponding blow up of the metric $g$ in a neighbourhood of $0$ in $B^+_1$, it follows from Proposition \ref{proposition foliation} that for all $i$ sufficiently large there exists a unique function $u_{i,t}=u(\hat{H}_i,t,g_i,w_i):D_{\theta}\rightarrow\real$ for each $|t|<\varepsilon_0$ such that $u_{i,t}=t+w_i$ on $\partial_+D_{\theta}$, the graph of $u_{i,t}$ over $D_{\theta}$ meets $\partial_0C_{\theta}$ orthogonally, it has constant mean curvature equal to $\hat{H}_i$ and its mean curvature vector points in the same direction as the mean curvature vector of $\hat\Sigma_i$.
Furthermore, $u_{i,t}$ varies smoothly on $t$ and defines a foliation of a region $D_{\theta}\times[-c,c]$ for some $c>0$ independent of $i$.

In case the mean curvature vector of $\hat\Sigma_i$ points in the negative direction of the $x_3$-axis, we use $u_{i,t}=u(-\hat{H}_i,t,g_i,w_i)$ instead.

Let $t_i\in(-\varepsilon_0,\varepsilon_0)$ be such that $u_{i,t_i}(0)=0$.
We claim that $\hat\Sigma_i$ is entirely below the graph of $u_{i,t_i}$ when $t_i\geq 0$ and reversely, $\hat\Sigma_i$ is entirely above the graph of $u_{i,t_i}$ when $t_i\leq 0$.

To prove it we use a standard application of maximum principle as follows {(see \cite{gilbarg-trudinger}*{Theorem 10.1} or \cite{m.taylorPDEIII}*{Chapter 14:Proposition 7.2} for a general version)}.
Since $\hat\Sigma_i$ converges graphically to $P_+$ away from $0$, we may find a function $v_i\in C^{2,\alpha}(D_\theta\setminus\{0\})\cap C^0(D_\theta)$ such that $\hat\Sigma_i$ is the graph of $v_i$ over $D_\theta$.
Denote by $g_{Eucl}$ the Euclidean metric on $C_\theta$ so that for any $X,Y\in T_{(x,t)C_{\theta}}$ we can write $g_i(X,Y)=g_{Eucl}(A_iX,Y)$ for some invertible matrix $A_i$ depending on $(x,t)$, where $(x,t)\in C_\theta$.
Given a function $f\in C^2(D_\theta\setminus\{0\})$ the normal vector of $\textup{graph}(f)=\{(x,f(x))\in C_\theta:x\in D_\theta$  with respect to $g_i$ is given by $\vec{N}_{\textup{graph}(f)}^{g_i}=\frac{1}{|A_i^{-1}(-\nabla f,1)|_{g_i}}A_i^{-1}(-\nabla f,1)$ so that $\Phi(f)=\textup{div}_{{D_\theta},g_i}\left(-\vec{N}_{\textup{graph}(f)}^{g_i}\right)=H_{\textup{graph}(f)}^{g_i}$ is the mean curvature of the graph of $f$ with respect to the metric $g_i$ and the linear elliptic operator
\begin{equation*}
L_if=\int_0^1\frac{d}{ds}|_{s=0}\Phi(\tau u_{i,t_i}+(1-\tau)v_i+sf)d\tau
\end{equation*}
can be written in divergence form $L_if=\textup{div}(X_i(x,\nabla f))$ where the function $x\mapsto X_i(x,\nabla f(x))$ is a vector field and $f\mapsto X_i(\cdot,\nabla f(\cdot))$ is linear.

Since the mean curvature vector of the graphs of $v_i$ and $u_{i,t_i}$ point in the same direction and have same magnitude, we have 
\begin{equation*}
\Phi(u_{i,t_i})-\Phi(v_i)=H^{g_i}_{\textup{graph}(u_{i,t_i})}-H^{g_i}_{\textup{graph}(v_{i})}=0.
\end{equation*}
Define $f_i=u_{i,t_i}-v_i$ so $f_i\in C^{2,\alpha}(D_\theta\setminus\{0\})\cap C^0(D_\theta)$ and
\begin{equation*}
\left\{
\begin{aligned}
Lf_i = & 0  , \text{ on } D_\theta\setminus\{0\}; \\
f_i = & t , \text{ on } \partial_+ D_\theta;\\
\frac{\partial f_i}{\partial x_1} = & 0 , \text{ on } \partial_0 D_\theta\setminus\{0\};\\
f_i(0) = & 0  .
\end{aligned}
\right.
\end{equation*}

Let us now consider the case $t\geq 0$ and the other case is similar.
Pick $x_0\in D_\theta$ such that $f_i(x_0)=\min f_i$ and suppose by contradiction that $f_i(x_0)<0$.
From the above, we have that $x_0\neq 0$ and since $t\geq 0$ then $x_0\not\in\partial_+D_\theta$.
Furthermore, it follows by the Hopf Lemma \cite{gilbarg-trudinger}*{Lemma 3.4} that $x_0\not\in\partial_0D_\theta\setminus\{0\}$.
Let us observe that $x_0$ is not on the corners of the domain either because those points also belong to $\partial_+D_\theta$.

In particular $x_0$ belongs to the interior of $D_\theta$ so, given an arbitrary domain $\Omega\subset D_\theta$ with smooth boundary containing $x_0$ we may apply maximum principle \cite{gilbarg-trudinger}*{Theorem 3.5} to get a contradiction thus proving that $f_i\geq 0$ on $D_\theta$.

Finally, pick another sequence $r_i'\rightarrow \infty$ and let $\hat\Sigma_i'$, $P_+'$ be its corresponding blow-up and tangent cone respectively.
For each $k$ we pick $i_k$ so that $r'_{i_k}>k r_k$.
Observe that $\hat\Sigma_{i_k}=\frac{r_{i_k}'}{r_k}\hat\Sigma_k$ which is contained to one side of the graph of $\frac{r_{i_k}'}{r_k}u_{k,t_k}$.
Since each $u_{k,t_k}$ is regular at $0$, $\frac{r_{i_k}'}{r_k}u_{k,t_k}$ must converge to a unique tangent cone, that is, $P_+$.
If $P_+'$ were different from $P_+$, then it would imply that $\frac{r_{i_k}'}{r_k}\hat\Sigma_k$ contains points on both sides of  the graph of $\frac{r_{i_k}'}{r_k}u_{k,t_k}$.
Notice that $P_+'$ must also be a half plane with free boundary on $\{x\in\real^3: x_1= 0\}$.
This concludes the claim.

From the above we know that the tangent cone of $\hat\Sigma$ at $0$ is unique and, by Lemma \ref{schauder estimates}, it is obtained as the limit of graphs of $C^{2,\alpha}$ functions converging in the $C^{2,\alpha}$ toplogy away from $0$.
In particular, for the fixed component $\hat\Sigma$ the tangent cone has multiplicity one and it is given by a half plane which we denote by $P_+$.
It follows that $\hat\Sigma\cup\{0\}$ is $C^1$ and $T^+_0(\hat\Sigma\cup\{0\})=P_+$.
Therefore, for $\epsilon>0$ sufficiently small we can find $u\in C^1(\bar D^+_\epsilon)\cap C^\infty(D^+_\epsilon\setminus \{0\})$, where $D^+_\epsilon\subset P_+$ is an open half-disk centered at $0$, so that $(\hat\Sigma\cup\{0\})\cap B^+_\epsilon$ can be written, under the appropriate identification, as the graph of $u$ in $T^+_0B_1^+$.
The function $u$ is a weak solution of $\divergent\left(\pm\frac{\nabla u}{\sqrt{1+|\nabla u|^2}}\right)=H_\Sigma$ in $D^+_\epsilon$ and $\frac{\partial u}{\partial \nu}=0$ on $\partial_0D^+_\epsilon$, where $\nu$ is the inward conormal of $\partial_0D^+_\epsilon$.

To prove regularity of $u$, we may use the fact that it satisfies a Neumann boundary condition to take an arbitrary prolongation $D_\epsilon$ of $D_\epsilon^+$ in $T_0B_1^+$ and extend $u$ by reflection to define a function $\hat u\in C^1(D_\epsilon)\cap C^\infty(D_\epsilon\setminus\partial_0D_\epsilon^+)$ that is a weak solution of $\divergent\left(\pm\frac{\nabla \hat u}{\sqrt{1+|\nabla \hat u|^2}}\right)=H_\Sigma$ in $D_\epsilon$.
Since we already know that $\hat u$ is $C^1$, it follows from interior elliptic regularity for quasilinear differential equations (see \cite{m.taylorPDEIII}*{Chapter 14:Theorem 4.4}) that $\hat u\in C^\infty(D_\epsilon)$ thus $u\in C^\infty(D^+_\epsilon)$.

To conclude the proof, we observe that the argument above holds for every connected component of $\Sigma\cap B_\delta^+$ that corresponds to a punctured half-disk.
That is, each component is smooth across $0$ and touching only at $0$.
Hence, $\Sigma\cup \{0\}$ is smooth properly almost embedded.
\end{proof}

We now state the equivalent result for removable interior singularities without repeating the proof.

\begin{theorem}\label{removable singularity interior}
Let $g$ be a Riemannian metric on $B_1$ and $\Sigma$ be a smooth, properly embedded, CMC surface in $B_1\setminus\{0\}$ with $H_\Sigma\leq H_0$, and $0\in\overline{\Sigma}$.
Suppose $\int_\Sigma|A^\Sigma|^2\leq C$ then $\Sigma\cup\{0\}$ is a smooth, properly almost embedded, CMC surface in $B_1$.
\end{theorem}

\begin{remark}
The proof is exactly the same with a minor modification on the construction of the foliation.
That is, by dropping the Neumann component on the definition of $\Phi$ in Proposition \ref{proposition foliation}.
\end{remark}


\section{Compactness Theorem}
In this section we will prove our main theorem for free boundary embedded CMC surfaces.
As we shall see, the limiting surface may not be embedded because CMC surfaces may have tangential self-intersection as long as the normal vector points at opposite directions.

\begin{theorem}\label{compactness fbcmc}
Let $N$ be a compact $3$-dimensional manifold with boundary.
Suppose $H_{\partial N}{>} H_0$ {with respect to the inward conormal of $N$ along $\partial N$} and let $\Sigma_i$ be a sequence of {connected} embedded free boundary CMC surfaces with mean curvature $H_i$, genus $g_i$ and number of ends $r_i$ satisfying:
\begin{enumerate}[(a)]
\item $|H_i|\leq H_0$;
\item $g_i\leq g_0$;
\item $r_i\leq r_0$;
\item $\textup{area}(\Sigma_i)\leq A_0$ and
\item $\textup{length}(\partial \Sigma_i)\leq L_0$.
\end{enumerate}

Then there exists a smooth properly almost embedded CMC surface $\Sigma\subset N$ and a finite set $\Gamma\subset \Sigma$ such that, up to a subsequence, $\Sigma_i$ converges to $\Sigma$ locally graphically in the $C^k$ topology on compact sets of $N\setminus \Gamma$ for all $k\geq 2$.
Moreover, if $\Sigma$ is minimal then it is properly embedded.

If in addition $\partial N$ satisfies $A^{\partial N}> 0$ with respect to the inward conormal of $\partial N$, then the convergence is $1$-sheeted away from $\Gamma$ whenever $H_\Sigma\neq 0$.
\end{theorem}
\begin{proof}
Let us denote by $A_i$ the second fundamental form of $\Sigma_i$.
Given $x\in \Sigma$, it follows from Gauss equation that $|A_i|^2(x) = H_i^2 + 2K_N(T_x\Sigma) - 2K_i(x)$, where $K_N, K_i$ are the sectional curvatures of $N$ and $\Sigma_i$ respectively.
From Gauss-Bonnet Theorem we have
\begin{equation*}
\int_{\Sigma_i}|A_i|^2 = H_i^2\area(\Sigma_i) + 2\int_{\Sigma_i}K_N(T_x\Sigma_i) +2\int_{\partial\Sigma_i}\kappa_g + 4\pi(2g_i+r_i -2),
\end{equation*} 
where $\kappa_g$ denotes the geodesic curvature of $\partial\Sigma_i$.
Because $\Sigma_i$ is free boundary, we have that that $\kappa_g = A_{\partial N}(\tau_{\partial\Sigma},\tau_{\partial\Sigma})$, where $A_{\partial N}$ is the second fundamental form of $\partial N$ with respect to the inward conormal vector of $N$ along $\partial N$ and $\tau_{\partial \Sigma}$ is the unit tangent vector of $\partial \Sigma$.
Since $N$ is compact, there exists a constant $C=C(N)>0$ such that
\begin{equation*}
\int_{\Sigma_i}|A_i|^2\leq C(H_i^2\area(\Sigma_i)+g_i+r_i+\area(\Sigma_i)+\textup{length}(\partial\Sigma_i)).
\end{equation*}
Hence, from hypotheses $(a)$-$(e)$ we have that the total curvature is uniformly bounded by a constant $C_0=C_0(N,H_0,g_0,r_0,A_0,L_0)>0$.

Denote by $\mu_i$ the Radon measure on $N$ defined by $\mu_i(U)=\int_{\Sigma_i\cap U}|A_i|^2\textup{dvol}_{\Sigma_i}$, for a subset $U\subset N$.
It follows from the above uniform bound that there exists a Radon measure $\mu$ in $N$ such that, up to a subsequence, $\mu_i$ converges weakly to $\mu$.
Furthermore, the set $\Gamma=\{p\in N: \mu(\{p\})\geq \epsilon_0\}$ has at most $C_0/\epsilon_0$ elements, here $\epsilon_0$ is taken from Theorem \ref{curvature estimates}.

For each $x\in N\setminus\Gamma$ there exists $r>0$ such that $\mu(\nball_r(x))<1$.
Hence, for each $i$ sufficiently large $\mu_i(\nball_r(x))<1$, that is, $\int_{\Sigma_i\cap\nball_r(x)}|A_i|^2<1$.
By possibly choosing a smaller value of $r$, it follows from Theorem \ref{curvature estimates} that
\begin{equation*}
\sup_{\Sigma_i\cap \nball_{\frac{r}{2}}}|A_i|^2\leq C,
\end{equation*}
for some constant $C>0$ independent of $i$.

Let $r<r_0$ as in Lemma \ref{schauder estimates} and suppose that, up to a subsequence, $\Sigma_i \cap \nball_{\frac{r}{4}}(x)$ is non-empty for all $i$ sufficiently large.
Since $\Sigma_i$ is embedded then $\Sigma_i \cap \nball_{\frac{r}{4}}(x)$ is the union of disjoint embedded connected components $\Sigma_{i,1},\ldots,\Sigma_{i,L}$ each of which is the graph of a function defined on an open ball of  fixed radius on $T_{y_{i,j}}\Sigma_{i,j}$ for some $y_{i,j}\in\Sigma_{i,j}, j=1,\ldots, L.$
Because $\Sigma_i$ is compact, the number of sheets $L$ must be finite and thus constant for $i$ sufficiently large.
Hence, under the appropriate identifications we may further assume that $\Sigma_{i,j}\cap\nball_{r'}(x)$ is the graph of a function $u_{i,j}$ defined on $B_{r'}(0)\subset T_{y_i}\Sigma_i$, for some $0<r'<\frac{r}{8}$ depending only on $L$ and a fixed $y_i\in\Sigma_i\cap\nball_{r'}(x)$.
Furthermore $u_{i,j}$ has uniform $C^{2,\alpha}$ bounds as in Lemma \ref{schauder estimates}.

We may assume that, up to a subsequence, $y_{i}$ converges to $y'$ and $T_{y_{i}}\Sigma_i$ converges to a plane $P\subset T_{y'}N$.
In which case, under further identifications, we have that $\Sigma_{i,j}\cap\nball_{\frac{r'}{2}}(y')$ is the graph of a function $u_{i,j}'$ defined on an open disk on $P$ (or half-disk in case $y'$ is on the boundary of $N$) and uniform $C^{2,\alpha}$ estimates for all $i$ sufficiently large, and each $j=1,\ldots, L$. 
Hence, up to a subsequence $u_{i,j}'$ converges to a function $u_j'$ in the $C^{2,\beta}$ topology for all $\beta<\alpha$ and $\Sigma_{i} \cap \nball_{\frac{r'}{2}}(y')$ converges to $\Sigma_{j}'=\graph( u_j')$ for each $j=1,\ldots, L$.
From which follows that $y'\in\Sigma'=\cup_j\Sigma_j'$ and $P$ is in fact $T_{y'}\Sigma'$.


Now, given any compact set $K\subset N\setminus \Gamma$ we cover it by finitely many open balls $\{\nball_{\frac{r}{4}}(x_k)\}_{k=0,\ldots,m}$ as above so that, up to a subsequence, $\Sigma_i\cap(\cup_k \nball_{\frac{r}{4}}(x_k))$ converges to a surface $\Sigma'\subset N\setminus\Gamma$ locally graphically in the $C^{2,\beta}$ topology on $K$.
Taking a countable exhaustion by compact sets and using a diagonal argument we have that, up to a subsequence, $\Sigma_i$ converges to $\Sigma'$ locally graphically in the $C^{2,\beta}$ topology on compact sets of $N\setminus\Gamma$.
Smooth convergence away from $\Gamma$ follows from Allard's regularity Theorem \cites{allard1972, allard1975, gruter-jost1986}.

Since $\Sigma_i$ is a properly embedded CMC surface with free boundary along $\partial N\setminus \Gamma$, then $\Sigma'$ is a properly almost embedded CMC surface in $N\setminus \Gamma$ with free boundary along $\partial N \setminus \Gamma$.
Observe that on a neighbourhood of any self-touching point of $\Sigma'$ the surface can be written as distinct embedded components that lie to one side of one another.
A transversal self-intersection is an open condition so it would contradict the fact that $\Sigma_i$ is embedded.
Furthermore, at a self-touching point $\Sigma'$ can be written only as two embedded components.
A third distinct embedded component would imply that the normal direction of at least two coincide thus contradicting the maximum principle.


Define $\Sigma$ to be the closure of $\Sigma'$, so $\Gamma=\Sigma\setminus\Sigma'$.
Since $\Gamma$ is finite, there exists $r > 0$ so that $\nball_r(p)\setminus\{p\}$ contains no points of $\Gamma$.
From graphical convergence on compact sets of $\nball_r(p)\setminus\{p\}$, each sheet must converge to an embedded component of $\Sigma\cap \nball_r(p)\setminus\{p\}$.
Since there are only finitely many sheets, we may pick $r>0$ sufficiently small so that every component of $\Sigma\cap \nball_r(p)\setminus\{p\}$ contains $p$ in its closure.

\begin{claim}
The limit surface $\Sigma$ is smooth and properly almost embedded.
\end{claim}
We know that $\int_{\Sigma\cap (\nball_r(p)\setminus\{p\})}|A_\Sigma|^2\leq C_0$.
Thus we may apply the removable singularity Theorems \ref{removable singularity boundary} or \ref{removable singularity interior} to each embedded component of $\Sigma\cap (\nball_r(p)\setminus\{p\})$ depending on whether $p$ belongs to the boundary of $N$ or to the interior.
Hence $\Sigma=\Sigma'\cup\Gamma$ is smooth and almost embedded.
Since $H_\Sigma\leq H_0 {<} H_{\partial N}$ it follows from the maximum principle that the interior of $\Sigma$ has no tangential touching points with $\partial N$, that is, $\Sigma$ is properly almost embedded.

Now we suppose that $H_\Sigma=0$,  then the maximum principle for minimal surfaces implies that there are no self-touching points and thus $\Sigma$ is embedded.

In the following we further assume that $A^{\partial N}>0$ with respect to the inward conormal vector of $N$ along $\partial N$ and prove the multiplicity part of the statement.

\begin{claim}\label{multiplicity claim}
If $H_\Sigma\neq 0$ then the convergence is $1$-sheeted away from $\Gamma$.
\end{claim}

If $\Gamma$ is empty then it follows from the fact that $\Sigma_i$ is connected and two-sided since $H_i>0$ for $i$ sufficiently large.

Let us consider only the boundary case and suppose $\Gamma\cap\partial \Sigma\neq\emptyset$.
The interior case follows the exact same argument.

First we will analyze the local picture around a point of singular convergence.
Take $p\in\Gamma\cap\partial\Sigma$ and $r>0$ sufficiently small so that $\nball_r(p)\cap\Gamma=\{p\}$ and $\nball_r(p)$ defines a geodesic ball adapted to the boundary and $\partial_+\nball_r(p)=\partial\nball_r(p)\cap(N\setminus\partial N)$ has mean curvature larger than $2H_\Sigma$.
Henceforth we shall omit the center of the ball in order to reduce notation.

It follows from Allard's regularity Theorem \cites{allard1972, allard1975, gruter-jost1986} that $p$ must be a self-touching point of $\Sigma$.
Therefore, by possibly taking $r>0$ smaller we may assume that $\Sigma\cap\nball_r=W_1\cup W_2$ where $W_m=W_m(p)$ is an embedded half disk meeting $\partial N$ orthogonally at $\partial_0W_m=\partial W_m\cap\partial N$ for each $m=1,2$.

Since the convergence is graphical away from $\Gamma$, for all $0<\delta<r$, $i$ sufficiently large and each $m=1,2$ there exist $L_m=L_m(p)$ positive integer, functions $u_{i,j}^m=u_{i,j}^{p,m}$, defined on $W_m\cap(\nball_r\setminus\nball_\delta)$ for $j=1,\ldots,L_m$ so that $\Sigma_i\cap V_\varepsilon^{X_m}(W_m\cap(\nball_r\setminus\nball_\delta))$ is given by the union of the graphs of $u_{i,j}^m$ for $j=1,\ldots,L_m$.

We define $\Lambda_{i,j}^m=\Lambda_{i,j}^m(p)$ to be the connected component of $\Sigma_i\cap\nball_r$ that contains the graph of $u_{i,j}^m$.

\begin{subclaim} If $\Lambda_{i,j}^1\cap(\nball_r\setminus\nball_\delta)$ is disconnected for some $j\in\{1,\ldots,L_1\}$, and for all $i$ sufficiently large then there exists $j'\in\{1,\ldots,L_2\}$ such that $\Lambda_{i,j}^1=\Lambda_{i,j'}^2$.
Equivalently, $\Lambda_{i,j}^1$ also contains the graph of $u_{i,j'}^2$ converging to $W_2\cap(\nball_r\setminus\nball_\delta)$.
\end{subclaim}
Since $\Lambda_{i,j}^1\cap(\nball_r\setminus\nball_\delta)$ is disconnected then $\Lambda_{i,j}^1$ must contain the graph of another function different from $u_{i,j}^1$.
Suppose the statement is false, that is, there exists $j'\in\{1,\ldots,L_1\}$ with $j'\neq j$ so that the graph of $u_{i,j'}^1$ is also contained in $\Lambda_{i,j}^1$.
Without loss of generality we may assume that $u_{i,j}^1>u_{i,j'}^1$.

Observe that $V_\varepsilon^{X_1}(W_1\cap(\nball_r\setminus\nball_\delta))\setminus(\graph(u_{i,j}^1)\cup\graph(u_{i,j'}^1))$ is given by $3$ connected components $V_{i,1},V_{i,2},V_{i,3}$ 
and $\nball_r\setminus\Lambda_{i,j}^1$ consists of only two connected components $Q_{i,1},Q_{i,2}$.
Since the mean curvature vector of $\Lambda_{i,j}^1$ never vanishes then it must either point into $Q_{i,1}$ everywhere or always point into $Q_{i,2}$. 

The neighbourhood $V_\varepsilon^{X_1}(W_1\cap(\nball_r\setminus\nball_\delta))$ is oriented by height and we may order the components so that $V_{i,1}$ lies above $\graph(u_{i,j}^1)$, $V_{i,2}$ lies between $\graph(u_{i,j}^1)$ and $\graph(u_{i,j'}^1)$, and $V_{i,3}$ lies below $\graph(u_{i,j'}^1)$.
Now, if the mean curvature vector of $\graph(u_{i,j}^1)$ points into $V_{i,1}$ then the mean curvature vector of $\graph(u_{i,j'}^1)$ must point into $V_{i,2}$ because both graphs converge smoothly to $W_1\cap(\nball_r\setminus\nball_\delta)$ so their mean curvature vectors must point in the same direction for $i$ sufficiently large.
Therefore, $V_{i,1}$ and $V_{i,2}$ must be contained in the same connected component of $\nball_r\setminus\Lambda_{i,j}^1$, either $Q_{i,1}$ or $Q_{i,2}$.
In particular we may construct a contractible closed path in $\nball_r$ with a single transversal intersection with $\Lambda_{i,j}^1$ which is a contradiction because $\Lambda_{i,j}^1$ is two-sided.

Similarly, if the mean curvature vector of $\graph(u_{i,j}^1)$ points into $V_{i,2}$ then the mean curvature vector of $\graph(u_{i,j'}^1)$ must point into $V_{i,3}$ which again leads to a contradiction.

\begin{subclaim}\label{one annulus}
For all $i$ sufficiently large there exists at most one connected component $\Sigma_{i,j_0}^1$ such that $\Lambda_{i,j_0}^1\cap\nball_r\setminus\nball_\delta$ is disconnected.
\end{subclaim}
Suppose by contradiction that there exist two such components, namely $\Lambda_{i,j_1}^1$ and $\Lambda_{i,j_2}^1$ with $j_1,j_2\in\{1,\ldots,L_1\}$ and $j_1\neq j_2$.

To simplify notation we will temporarily omit the sequence lower index $i$ whenever it is not needed and write $\Lambda_m=\Lambda_{i,j_m}^1$ for each $m=1,2$.

%

\begin{figure}[H]
\includegraphics[scale=0.5]{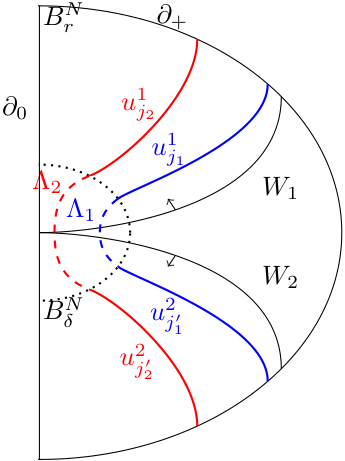}
\end{figure}
 
\justifying\parindent=12pt

To each $\Lambda_m$, $m=1,2$, we have a pair of functions $u_{j_m}^m,u_{j_m'}^m$ with $j_m\in\{1,\ldots,L_1\}$ and $j_m'\in\{1,\ldots,L_2\}$ such that $u_{j_m}^1,u_{j_m}^1$ converge to $W_1\cap(\nball_r\setminus\nball_\delta)$ and $u_{j_m'}^2,u_{j_m'}^2$ converge to $W_2\cap(\nball_r\setminus\nball_\delta)$.
Observe that the graphs are taken with respect to an extension of the normal vectors of each $W_1$ and $W_2$ so we may further assume that $u_{j_1}^1<u_{j_2}^1$ and $u_{j_1'}^2<u_{j_2'}^2$

Denote by $U$ the connected component of $\nball_r\setminus(\Lambda_{1}\cup\Lambda_{2})$ that lies between the two surfaces.
By the above assumption, the mean curvature of $\Lambda_{1}$ points into $U$ and the mean curvature of $\Lambda_{2}$ points away from $U$.
Define $\mathfrak{F}$ to be the family of subsets $\Omega\subset U$ such that $\Omega$ has finite perimeter, $\partial\Omega$ is a rectifiable $2$-current and $\Lambda_1\subset\partial\Omega$.
For each $\Omega\in\mathfrak{F}$ we denote $\partial_0\Omega=\partial\Omega\cap\partial N$ and $S=\partial\Omega\setminus(\partial_0\Omega\cup\Lambda_1)$.
In particular $\partial_+S=\partial S\cap\partial_+\nball_r=\partial_+\Lambda_1$ and $\partial(\partial_0\Omega)=\partial_0 S\cup\partial\Lambda_1$.
That is, $S$ and $\Lambda_1$ are homologous in $U$ relative to $\partial_0 U=\partial U\cap\partial N$.

We define the functional $F_{i}$ defined on $\mathfrak{F}$ as
\begin{equation*}
F_i(\Omega)=\area(S)+H_i\volume(\Omega).
\end{equation*}

The goal is to find a minimiser for $F_i$ strictly contained in $U$.
With that in mind we first construct barriers along $\partial U\setminus\partial_0 U$ to ensure that the minimiser is in the interior of $U$.

Firstly, for each $m=1,2$ we choose $\vec{N}_m=\frac{\vec{H}_{\Lambda_m}}{|\vec{H}_{\Lambda_m}|}$ as the normal vectorfield of $\Lambda_m$.
Let $Z_m$ be an arbitrary vectorfield extension of $\vec{N}_m$ near $\Lambda_m$torfield with compact support such that $Z_m|_{\partial N}\in T\partial N$ for each $m=1,2$ and $\support(Z_1)\cap\support(Z_2)=\emptyset$.

We observe that neither $\Lambda_1$ or $\Lambda_2$ are strongly stable as free boundary CMC surfaces with respect to the area functional.
In fact, if $\Lambda_m$ were strongly stable for either $m=1,2$ then it would follow from \cite{hong-saturnino2023}*{Theorem 1.6} (see Theorem \ref{curvature estimates 2}) that $\Lambda^1_{i,j_m}$ has uniform curvature estimates, independent of the sequence index $i$, which we briefly reintroduce to illustrate this step of the argument.
Therefore, it converges, up to a subsequence, graphically smoothly to $W_1\cup W_2$ as $i$ tends to infinity.
It follows from Claim 2.1 that $\partial_+\Lambda^1_{i,j_m}$ converges smoothly to $\partial_+W_1\cup \partial_+W_2$, which implies that $\Lambda^1_{i,j_m}$ must converge graphically smoothly to $W_1\cup W_2$ everywhere including at $p$.
Since the normal vectors of $W_1$ and $W_2$ point at opposite directions on $p$, it implies that $\Lambda^1_{i,j_m}$ must be given as the graph of a function over $W_1$ and another function over $W_2$, but $\Lambda^1_{i,j_m}$ is connected, which implies that these functions must coincide at some point of $W_1\cap\nball_\delta$ and $W_2\cap\nball_\delta$, which contradicts the fact that $\Sigma_i$ is embedded.

Therefore, for each $m=1,2$, we may find positive eigenfunctions $\varphi_m$ solving
\begin{equation*}
\left\{
\begin{aligned}
&\Delta_{\Lambda_m}\varphi_m + (|A_{m}|^2+Ric^N(\vec{N}_{m},\vec{N}_{m})+\lambda_{1,m})\varphi_m=0 \textup{ on }\Lambda_m,\\
&\varphi_m = 0  \textup{ on } \partial_+\Lambda_m \textup{ and }\\
&\frac{\partial}{\partial\eta_m}\varphi_m=0 \textup{ on } \partial_0\Lambda_m,
\end{aligned}
\right.
\end{equation*}
where $A_{m}$ is the second fundamental form of $\Lambda_m$ with respect to $\vec{N}_m$, $\lambda_{1,m}<0$ is the first eigenvalue of the Jacobi operator of $\Lambda_m$ in $\nball_r$ with Neumann conditions on $\partial_0\Lambda_m$ and $\eta_m$ is the inward conormal of $\Lambda_m$ along $\partial_0 \Lambda_m$.
Since $\varphi_m$ is the first eigenfunction we further have that $\varphi_m>0$ on $\Lambda_m\setminus\partial\Lambda_m$.

By possible perturbing $r>0$ we may further assume that $0$ is not an eigenvalue of the Jacobi operator of $\Lambda_m$ in $\nball_n$ for both $m=1,2$.
Hence we can find a function $\xi_m$ satisfying
\begin{equation*}
\left\{
\begin{aligned}
&\Delta_{\Lambda_m}\xi_m + (|A_m|^2+Ric^N(\vec{N}_{m},\vec{N}_{m})\xi_m=1 \textup{ on }\Lambda_{m},\\
&\xi_m = 0  \textup{ on } \partial_+\Lambda_m \textup{ and }\\
&\frac{\partial}{\partial\eta_m}\xi_m=0 \textup{ on } \partial_0\Lambda_m.
\end{aligned}
\right.
\end{equation*}
It follows from the Hopf Lemma \cite{gilbarg-trudinger}*{Lemma 3.4} that $\frac{\partial}{\partial\eta^{\partial_+\Lambda_m}_{in}}\varphi_m>0$, where $\eta^{\partial_+\Lambda_m}_{in}$ is the inward conormal of $\partial_+\Lambda_m$.
Therefore, we may find $c>0$ depending on $\xi_m$ and normal derivatives of $\varphi_m$ along $\partial_+\Lambda_m$ such that $v_m=\varphi_m+c\xi_m$ is positive on $\Lambda_m\setminus\partial_+\Lambda_m$.

Now we write $\Lambda_m(t)=\{\Phi^{Z_m}(q,tv_m(q)):q\in\Lambda_m\}$.
For $|t|$ sufficiently small $\Lambda_m(t)$ defines a foliation of an open region in $N$.
Furthermore, when $m=1$ we may find $\rho_1>0$ sufficiently small such that $U_1=\{\Phi^{Z_1}(q,tv_1(q)):q\in\Lambda_1,0<t\leq\rho_1\}\subset U$ and when $m=2$ we may find $\rho_2>0$ sufficiently small such that $U_2=\{\Phi^{Z_2}(q,tv_2(q)):q\in\Lambda_2,-\rho_2\leq t<0\}\subset U$.

We write $H_m(t)$ for the mean curvature of $\Lambda_m(t)$ and use the linearization of the mean curvature to obtain
\begin{equation*}
\begin{aligned}
\frac{d}{dt}|_{t=0}H_m(t) & =\Delta_{\Lambda_m}v_m + (|A_m|^2+Ric^N(\vec{N}_{m},\vec{N}_{m}))v_m\\
                          & = -\lambda_{1,m}\varphi_m+c > 0
\end{aligned}
\end{equation*}

Let $Y_m$ be the vectorfield on $\{\Phi^{Z_m}(q,tv_m(q)):q\in\Sigma_m,-\rho_m\leq t\leq\rho_m\}$ so that $Y_m(\Phi^{Z_m}(q,tv_m(q)))$ is the normal vector to $\Lambda_m(t)$ at $\Phi^{Z_m}(q,tv_m(q))$ and $Y_m(q)=\vec{N}_m(q)$ at $t=0$.
In particular, $H_m(t)=-\divergent_N(Y_m)$ and since $H_m(0)=H_i$, by possibly making $\rho_m>0$ smaller we may assume that 
\begin{equation*}\tag{*}
\begin{aligned}
\divergent_N(Y_1) & < -H_i \text{ on } U_1, \text{ that is, }0<t\leq\rho_1 \text{ and }\\
\divergent_N(Y_2) & > -H_i \text{ on } U_2, \text{ that is, }-\rho_2\leq t< 0.
\end{aligned}
\end{equation*}
Note that $\Lambda_m(t)$ does not necessarily intersect $\partial N$ perpendicularly, however we compute
$\frac{d}{dt}|_{t=0}g(Y_m,\nu^{\partial N}_{in})=-v_mA^{\partial N}(\vec{N}_{m},\vec{N}_{m})<0$,
where $\nu^{\partial N}_{in}$ is the inward conormal vector of $N$ along $\partial N$.
Since $g(Y_m,\nu^{\partial N}_{in})=0$ at $t=0$, by possibly making $\rho_m>0$ smaller, we may further assume that
\begin{equation*}\tag{**}
\begin{aligned}
g(Y_1,\nu^{\partial N}_{in}) & < 0 \text{ on } U_1, \text{ that is, }0<t\leq\rho_1 \text{ and }\\
g(Y_2,\nu^{\partial N}_{in}) & > 0 \text{ on } U_2, \text{ that is, }-\rho_2\leq t< 0.
\end{aligned}
\end{equation*}
That concludes the construction of the barrier around $\Lambda_1$ with $U_1$ and $S_1=\Lambda_1(\rho_1)$ and the barrier around $\Lambda_2$ with $U_2$ and $S_2=\Lambda_2(-\rho_2)$.
Since $Y_m$ coincides with the normal vectors of $\Sigma_m$ at $t=0$, we may further assume, by possibly making $\rho_m>0$ smaller, that $Y_1$ points away from $U_1$ along $S_1$ and $Y_2$ points towards $U_2$ along $S_2$.

Finally we note that $\nball_r\setminus\Lambda_1$ has two connected components and we choose $\tilde{U}$ to be the connected component with mean convex boundary, that is, $\vec{H}_1$ points into $\tilde{U}$.
Furthermore, $\partial\tilde{U}$ has three components, $\partial_0\tilde{U}=\tilde{U}\cap\partial N$, $\Lambda_1$ and $\partial_+\nball_r\cap\tilde{U}$ all of which have mean curvature vector pointing into $\tilde{U}$ and each of which have scalar mean curvature larger than $\frac{H_\Sigma}{2}>0$ as long as $i$ is sufficiently large.
It follows that we can find an area minimising surface $S_0$ with boundary such that $S_0$ is properly embedded in $\tilde{U}$ (see for example \cite{white2010}), $S_0$ intersects $\partial_0\tilde{U}$ perpendicularly and $\partial_+ S_0=\partial_+\Lambda_1$.
Let $U_0$ be the connected component of $\tilde{U}\setminus(\Lambda_1\cup S_0)$ bounded by $S_0$, $\Sigma_1$ and $\partial_0 \tilde{U}$.


%

\begin{figure}[H]
\includegraphics[scale=0.5]{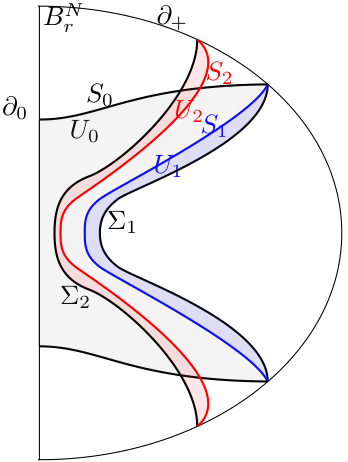}
\end{figure}


\justifying\parindent=12pt

With the barriers constructed above (see the diagram above for an illustrative picture of these barriers), we will prove that it is always possible to find a minimiser of $F_i$ such that the corresponding partial boundary surface lies in the region $(U\setminus(U_1\cup U_2))\cap U_0$.

\begin{subclaim}\label{barrier inequalities}
If $\Omega\in\mathfrak{F}$ has $S=\partial\Omega\setminus(\partial_0\Omega\cup\Lambda_1)$ smooth and transversal to $S_0$, $S_1$ and $S_2$ whenever the respective intersection is non-empty, then the following hold:
\begin{enumerate}[(a)]
\item If $\Omega\not\subset U_{0}$ then $F_i(\Omega\cap U_{0})<F_i(\Omega)$;
\item If $U_{1}\not\subset\Omega$ then $F_i(\Omega\cup U_{1})<F_i(\Omega)$;
\item If $\Omega\cap U_{2}\neq\emptyset$ then $F_i(\Omega\setminus U_{2})<F_i(\Omega)$;
\end{enumerate}
\end{subclaim}

%
%
%

\begin{figure}[H]
\includegraphics[scale=0.5]{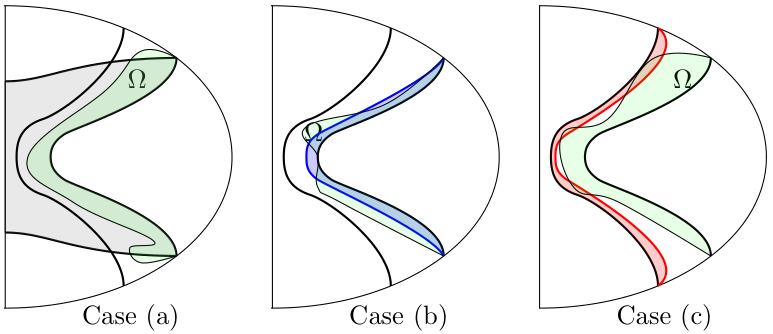}
\end{figure}

\justifying\parindent=12pt

\noindent\textbf{Case (a).} Let $\Omega'=\Omega\cap U_0$ so that $S'=\partial\Omega'\setminus(\partial_0\Omega'\cup\Lambda_1)$ can be written as a disjoint union $S'=(S_0\cap\Omega)\cup(S\cap U_0)$.
Since $S_0$ is area minimising and $\partial(S_0\cap\Omega)=\partial(S\setminus U_0)$ we have that $\area(S_0\cap\Omega)\leq\area(S\setminus U_0)$ and $\Omega'\subset\Omega$ is, by assumption, a strict inclusion, so that $\vol(\Omega')<\vol(\Omega)$.
Hence,
\begin{equation*}
\begin{aligned}
F_i(\Omega') & = \area(S')+H_i\vol(\Omega')\\
             & = \area(S_0\cap\Omega)+\area(S\cap U_0)+H_i\vol(\Omega')\\
						 & < \area(S\setminus U_0)+\area(S\cap U_0)+H_i\vol(\Omega)\\
						 & = \area(S)+H_i\vol(\Omega)\\
						 & = F_i(\Omega).
\end{aligned}
\end{equation*}

\noindent\textbf{Case (b).} Let $\Omega'=\Omega\cup U_1$ so that $\Omega'=\Omega\cup(U_1\setminus\Omega)$ and $S'=\partial\Omega'\setminus(\partial_0\Omega'\cup\Lambda_1)$ can be written as $S'=(S\setminus U_1)\cup(S_1\setminus\Omega)$.
Recall from $(*)$ in the construction of the barriers that $\divergent_N(Y_1)<-H_i$ on $U_1$. 
If we denote by $\eta_{out}$ the outer normal vector along $\partial(U_1\setminus\Omega)$, we note that $\partial(U_1\setminus\Omega)=\partial_0(U_1\setminus\Omega)\cup(S_1\setminus\Omega)\cup(S\cap U_1)$ and use the divergence theorem to compute:
\begin{equation*}
\begin{aligned}
-H_i\vol(U_1\setminus\Omega) & > \int_{U_1\setminus\Omega}\divergent_N(Y_1)\\
                             & = \int_{\partial_0(U_1\setminus\Omega)}g(Y_1,\eta_{out})+\int_{S_1\setminus\Omega)}g(Y_1,\eta_{out})+\int_{S\cap U_1}g(Y_1,\eta_{out}).
\end{aligned}
\end{equation*}
We observe that $\eta_{out}=-\nu^{\partial N}_{in}$ along $\partial_0(U_1\setminus\Omega)$ and by $(**)$ we have $g(Y_1,-\nu^{\partial N}_{in})>0$ on $\partial_0(U_1\setminus\Omega)$.
Since $Y_1$ points away from $U_1$ on $S_1$, it follows that $\eta_{out}=Y_1$ along $S_1\setminus\Omega$, so $g(Y_1,\eta_{out})=1$ on $S_1\setminus\Omega$.
Furthermore, we have trivially that $g(Y_1,\eta_{out})\geq -1$ on $S\cap U_1$.
Thus,
\begin{equation*}
-H_i\vol(U_1\setminus\Omega) > 0+\area(S_1\setminus\Omega)-\area(S\cap U_1).
\end{equation*}
Which we rearrange to obtain:
\begin{equation*}
H_i\vol(U_1\setminus\Omega) <-\area(S_1\setminus\Omega)+\area(S\cap U_1).
\end{equation*}
We conclude by computing:
\begin{equation*}
\begin{aligned}
F_i(\Omega') & = \area(S') + H_i\vol(\Omega')\\
             & = \area(S\setminus U_1)+\area(S_1\setminus\Omega)\\
						 & \quad + H_i\vol(U_1\setminus\Omega)+ H_i\vol(\Omega)\\
						 & < \area(S\setminus U_1)+\area(S_1\setminus\Omega)\\
             & \quad -\area(S_1\setminus\Omega)+\area(S\cap U_1)+ H_i\vol(\Omega)\\
						 & = \area(S)+H_i\vol(\Omega)\\
						 & = F_i(\Omega).
\end{aligned}
\end{equation*}

\noindent\textbf{Case (c).} Let $\Omega'=\Omega\setminus U_2$ so that $\Omega=\Omega'\cup(\Omega\cap U_2)$ and $S'=\partial\Omega'\setminus(\partial_0\Omega'\cup\Lambda_1)$ can be written as $S'=(S\setminus U_2)\cup(S_2\cap\Omega)$.
Recall from $(*)$ in the construction of the barriers that $\divergent_N(Y_2)>-H_i$ on $U_2$.
If we denote by $\eta_{out}$ the outer normal vector along $\partial(U_2\cap\Omega)$, we note that $\partial(U_2\cap\Omega)=\partial_0(U_2\cap\Omega)\cup(S_2\cap\Omega)\cup(S\cap U_2)$ and use the divergence theorem to compute:
\begin{equation*}
\begin{aligned}
-H_i\vol(U_2\cap\Omega) & < \int_{U_2\cap\Omega}\divergent_N(Y_2)\\
                             & = \int_{\partial_0(U_2\cap\Omega)}g(Y_2,\eta_{out})+\int_{S_2\cap\Omega)}g(Y_2,\eta_{out})+\int_{S\cap U_2}g(Y_2,\eta_{out}).
\end{aligned}
\end{equation*}
We observe that $\eta_{out}=-\nu^{\partial N}_{in}$ along $\partial_0(U_2\cap\Omega)$ and by $(**)$ we have $g(Y_2,-\nu^{\partial N}_{in})<0$ on $\partial_0(U_2\cap\Omega)$.
Since $Y_2$ points into $U_2$ on $S_2$, it follows that $\eta_{out}=-Y_2$ along $S_2\cap\Omega$, so $g(Y_2,\eta_{out})=-1$ on $S_2\cap\Omega$.
Furthermore, we have trivially that $g(Y_2,\eta_{out})\leq 1$ on $S\cap U_2$.
Therefore,
\begin{equation*}
-H_i\vol(U_2\cap\Omega) < 0 - \area(S_2\cap\Omega)+\area(S\cap U_2).
\end{equation*}
Which we rearrange to obtain:
\begin{equation*}
\area(S_2\cap\Omega) < \area(S\cap U_2)+H_i\vol(U_2\cap\Omega).
\end{equation*}
We conclude by computing:
\begin{equation*}
\begin{aligned}
F_i(\Omega') & = \area(S') + H_i\vol(\Omega')\\
             & = \area(S\setminus U_2)+\area(S_2\cap\Omega)\\
						 & \quad + H_i\vol(\Omega')\\
						 & < \area(S\setminus U_2)+\area(S\cap\Omega)+H_i\vol(U_2\cap\Omega)\\
             & \quad +H_i\vol(\Omega')\\
						 & = \area(S)+H_i\vol(\Omega)\\
						 & = F_i(\Omega).
\end{aligned}
\end{equation*}
Which concludes the proof of Claim \ref{barrier inequalities}.

We now conclude the proof of Claim \ref{one annulus}.
Let us reintroduce the sequence index $i$ and recall $\Lambda_m=\Lambda_{i,j_m}^1$ for each $m=1,2$.

Take a minimising sequence $\{\Omega_{i,l}\}_{l\in\naturals}$ of $F_i$ for which we may assume that $\Omega_{i,l}\subset U_{0}$, $U_1\subset\Omega_{i,l}$ and $\Omega_{i,l}\cap U_2 = \emptyset$ for all $l$ sufficiently large.
It follows from compactness and regularity (see \cite{cook1985}) that $\Omega_{i,l}$ converges, up to a subsequence, to a regular $\Omega_{i,\infty}$ so that $S_i=\partial\Omega_{i,\infty}\setminus(\partial_0\Omega_{i,\infty}\cup\Lambda_{i,j_1}^1)$ is a regular weakly stable properly embedded CMC surface in $U$ with free boundary in $\partial_0 U$ and homologous to $\Lambda_{i,j_1}^1$ relative to $\partial_0 U$.
Since $S_i$ is weakly stable, it follows from \cite{hong-saturnino2023}*{Theorem 1.6} (see Theorem \ref{curvature estimates 2}) that it has uniform curvature estimates.
We can also compute $\area(S_i)\leq F_i(\Omega_{i,\infty})\leq F_i(\Lambda_{i,j_1}^1)=\area(\Lambda_{i,j_1}^1)\leq 2\area(W_1\cup W_2)$, so $S_i$ also has uniform area bounds.
Therefore, it must converge, up to a subsequence, graphically smoothly to $W_1\cup W_2$ since $S_i$ is homologous to $\Lambda_{i,j_1}^1$ and, by Claim 2.1, $\partial_+S_i=\partial_+\Lambda_{i,j_1}^1$ converges smoothly to $\partial_+W_1\cup \partial_+W_2$.
However, arguing similarly to the proof of Claim 2.1, it would imply that either $\Lambda_{i,j_2}^1$ converges graphically everywhere or that $S_i$ intersect $\Lambda_{i,j_2}^1$, which are both contradictions.
This concludes the proof of Claim \ref{one annulus}.

Finally we will finish the proof of Claim \ref{multiplicity claim}.
First we write $\Sigma$ as the (non-disjoint) union of embedded connected components $\{\Sigma^{(k)}\}_{k=1,\ldots,l}$ for some positive integer $l$ such that each pair of embedded components may only intersect on the self-touching set $\mathcal{S}$ of $\Sigma$.
Let $r>0$ be sufficiently small such that $\nball_r(p)\cap\Gamma=\{p\}$ for all $p\in\Gamma$ and $0<\delta<r$ as in the previous claims.

Now, for each $k=1,\ldots,l$ we may take an extension of the normal vectorfield given by the direction of the mean curvature vector $\vec{H}^{(k)}$ of $\Sigma^{(k)}$ appropriately adapted to $\partial N$ and a neighbourhood $V^{(k)}$ of $\Sigma^{(k)}$ as in Definition \ref{definition convergence} such that $\Sigma_i\cap( V^{(k)}\setminus(\cup_{p\in\Gamma}\nball_\delta(p)))$ can be written as the graph of functions $u^{(k)}_{i,1}>\ldots>u^{(k)}_{i,m_k}$ for some positive integer $m_k$.
To each function $u^{(k)}_{i,j}$ we associate a connected component $\Sigma^{(k)}_{i,j}$ of $\Sigma_i\cap V^{(k)}\nball_\delta(p)))$ that contains its graph.

We are assuming by contradiction that $m_k\geq 2$ for some fixed $k$, in which case we may also infer that $\Sigma^{(k)}\cap\Gamma\neq\emptyset$, otherwise the convergence would be smooth everywhere on $\Sigma^{(k)}$ and thus $\Sigma_i\cap V^{(k)}$ would contain at least two graphical sheets and thus be disconnected, which is a contradiction.
Let us fix the first sheet ordered by height $\Sigma^{(k)}_{i,1}$ and for each $p\in\Sigma^{(k)}\cap\Gamma$ let $W_1(p),W_2(p)$ be the corresponding disks or half-disks (depending on whether $p$ belongs to the boundary or not).
By possibly relabelling $W_1(p),W_2(p)$ we may further assume that the embedded component $\Sigma^{(k)}$ contains $W_1(p)$ for all $p\in\Sigma^{(k)}\cap\Gamma$.
As in the beginning of Claim 2 we take $\Lambda^{(k)}_{i,1}(p)$ to be the connected component of $\Sigma_i\cap\nball_r(p)$ that contains $\Sigma^{(k)}_{i,1}\cap(\nball_r(p)\setminus\nball_\delta(p))$.
See the diagrams below for two possibilities of the above description.

%
\begin{figure}[H]
\includegraphics[scale=0.5]{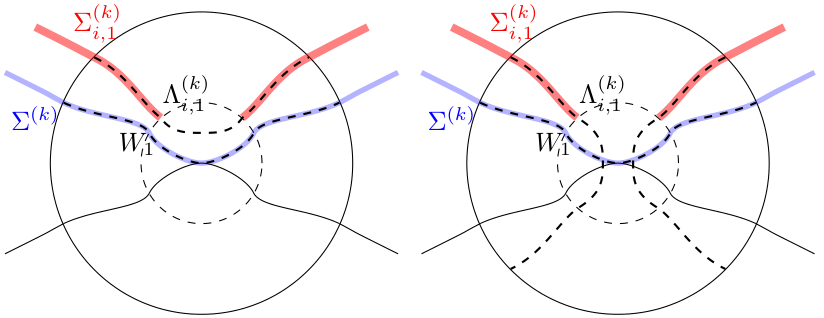}
\end{figure}
 
\justifying\parindent=12pt

\begin{subclaim}
There exists at least one $p\in\Sigma^{(k)}\cap\Gamma$ such that $\Lambda^{(k)}_{i,1}(p)\cap(\nball_r(p)\setminus\nball_\delta(p))$ is disconnected.
\end{subclaim}

Suppose that for every $p\in\Sigma^{(k)}\cap\Gamma$ we have that $\Lambda^{(k)}_{i,1}(p)\cap(\nball_r(p)\setminus\nball_\delta(p))$ is connected, so that $\Lambda^{(k)}_{i,1}(p)\cap(\nball_r(p)\setminus\nball_\delta(p))$ converges to $W_1(p)\cap(\nball_r(p)\setminus\nball_\delta(p))$ for every $p\in\Sigma^{(k)}\cap\Gamma$.
In this case $\Sigma^{(k)}_{i,1}\cup(\cup_{p\in\Sigma^{(k)}\cap\Gamma}\Lambda^{(k)}_{i,1}(p))$ would correspond to a (global) connected component of $\Sigma_i$ hence disjoint from the other sheets, which contradicts the fact that $\Sigma_i$ is connected.

%

\begin{figure}[H]
\includegraphics[scale=0.5]{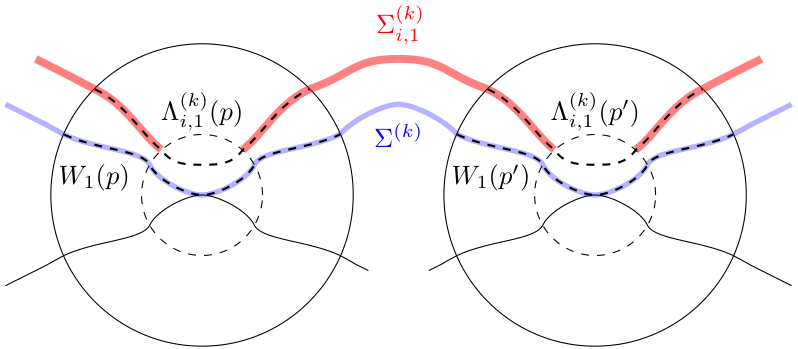}
\end{figure}
 
\justifying\parindent=12pt

Now we fix the point $p$ so that $\Lambda^{(k)}_{i,1}(p)\cap(\nball_r(p)\setminus\nball_\delta(p))$ is disconnected and will omit it from notation.
It follows from Claim 2.1 that we may label the connected components of $\Lambda^{(k)}_{i,1}\cap(\nball_r\setminus\nball_\delta)$ as $\Lambda^{1,(k)}_{i,1}$ and $\Lambda^{2,(k)}_{i,1}$ so that $\Lambda^{1,(k)}_{i,1}$ converges to $W_1\cap(\nball_r\setminus\nball_\delta)$ and $\Lambda^{2,(k)}_{i,1}$ converges to $W_2\cap(\nball_r\setminus\nball_\delta)$.

From Claim \ref{one annulus} we know that all other components of $\Sigma_i\cap\nball_r$ are connected on $\nball_r\setminus\nball_\delta$.
Take the next lower sheet $\Sigma^{(k)}_{i,2}$ (ordered by height on $V^{(k)}\setminus(\cup_{p\in\Gamma}\nball_\delta(p))$) and $\Lambda^{(k)}_{i,2}$ the connected component of $\Sigma_i\cap\nball_r$ that contains $\Sigma^{(k)}_{i,2}\cap(\nball_r\setminus\nball_\delta)$.
Then $\Lambda^{(k)}_{i,2}\cap(\nball_r\setminus\nball_\delta)$ must be connected and  converges to $W_1\cap(\nball_r\setminus\nball_\delta)$.
However, $\Lambda^{(k)}_{i,2}\cap(\nball_r\setminus\nball_\delta)$ must lie below $\Lambda^{1,(k)}_{i,1}$ but it lies above $\Lambda^{(2,k)}_{i,1}$ because the latter converges to $W_2\cap(\nball_r\setminus\nball_\delta)$.
Therefore we may pick two points, one in $\Lambda^{1,(k)}_{i,1})$ and another in $\Lambda^{2,(k)}_{i,1})$ and join them by a path that intersects $\Lambda^{(k)}_{i,2})$ transversally at a single point.
We may further connect those points by another path entirely contained in $\Lambda^{(k)}_{i,1}$ hence not intersecting $\Lambda^{(k)}_{i,2}$.
The concatenation of the two paths above gives us a closed contractrible path in $\nball_r$ that intersects $\Lambda^{(k)}_{i,2}$ transversally at a single point which is a contradiction since $\Lambda^{(k)}_{i,2}$ separates $\nball_r$ into two disjoint connected components.

%

\begin{figure}[H]
\includegraphics[scale=0.5]{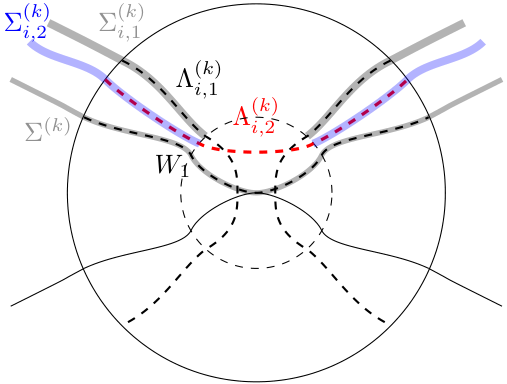}
\end{figure}
 
\justifying\parindent=12pt
Since the initial choice of embedded connected component of $\Sigma$ was arbitrary we conclude that the convergence must be one sheeted away from $\Gamma$ for every connected embedded component, which concludes the proof of Claim 2 and finishes the theorem.

\end{proof}

\begin{remark}\label{mistake 1 remark}
In the closed case (i.e, when neither the ambient space or the sequence has boundary) the proof of Claim 2 can be repeated ipsis litteris for hypersurfaces when the ambient dimension is $3\leq n \leq 7$.
The only steps we use dimension $2$ in this part of the argument was in the regularity of free boundary CMC surfaces that minimize $F_i$ and the curvature estimates for weakly stable free boundary CMC surfaces.
The closed case has equivalent versions of these statements valid for hypersurface dimensions $2\leq n-1\leq 6$ (see for example \cite{bellettini-chodosh-wickramasekera2019,schoen-simon-almgren1977}).
That said, it corrects the minor mistake found in the Multiplicity Analysis proof of \cite{bourni-sharp-tinaglia}*{Theorem 4.1, Claims 4.3-4.5}.
Our proof is entirely adapted from their original paper but we make a small change in the construction of the barriers.
\end{remark}

\begin{remark}\label{mistake 2 remark}
When the limit surface is minimal it is not yet clear that the number of sheets must be limited.
It is possible that near a point of singular convergence there are multiple components that are disconnected after removing a small ball, that is, Claim 2.2 may fail.
See the diagram below for a diagram that our proof is unable to exclude.

%

\begin{figure}[H]
\includegraphics[scale=0.5]{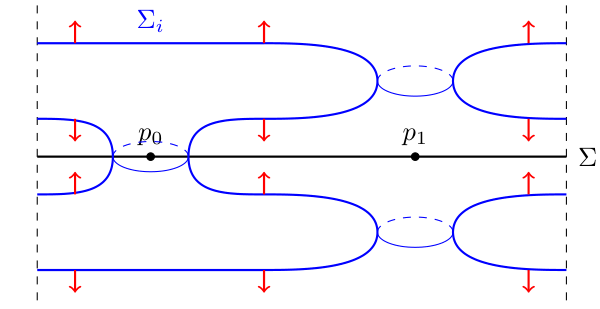}
\end{figure}
 
\justifying\parindent=12pt

The local picture around $p_1$ enables multiple sheets even when the sequence is connected but it can be excluded when the mean curvature of $\Sigma$ is positive thanks to Claim 2.1.
\end{remark}

\bibliographystyle{plain}
\bibliography{bibliography}

\end{document}